\documentclass[twocolumn,english,journal]{IEEEtran}
\usepackage[T1]{fontenc}
\usepackage{babel}
\usepackage{array}
\usepackage{float}
\usepackage{multirow}
\usepackage{amsmath}
\usepackage{amsthm}
\usepackage{amssymb}
\usepackage{graphicx}
\usepackage[numbers]{natbib}
\usepackage[unicode=true,
 bookmarks=true,bookmarksnumbered=true,bookmarksopen=true,bookmarksopenlevel=1,
 breaklinks=false,pdfborder={0 0 0},pdfborderstyle={},backref=false,colorlinks=false]
 {hyperref}
\hypersetup{pdftitle={Your Title},
 pdfauthor={Your Name},
 pdfpagelayout=OneColumn, pdfnewwindow=true, pdfstartview=XYZ, plainpages=false}

\makeatletter

\providecommand{\tabularnewline}{\\}
\floatstyle{ruled}
\newfloat{algorithm}{tbp}{loa}
\providecommand{\algorithmname}{Algorithm}
\floatname{algorithm}{\protect\algorithmname}

\theoremstyle{plain}
\newtheorem{thm}{\protect\theoremname}
\theoremstyle{plain}
\newtheorem{lem}[thm]{\protect\lemmaname}

\usepackage{algorithmic}
\usepackage[bb=boondox]{mathalfa}
\makeatletter
\def\thmhead@plain#1#2#3{%
  \thmname{#1}\thmnumber{\@ifnotempty{#1}{ }\@upn{#2}}%
  \thmnote{ {\the\thm@notefont#3}}}
\let\thmhead\thmhead@plain
\makeatother
\usepackage{xcolor}

\usepackage{subfig}
\usepackage{babel}
\makeatletter

\usepackage{array}
\usepackage{boldline}
\usepackage{multirow}
\usepackage{nicematrix}
\usepackage{threeparttable}

\@ifundefined{showcaptionsetup}{}{%
 \PassOptionsToPackage{caption=false}{subfig}}
\usepackage{subfig}
\makeatother

\providecommand{\lemmaname}{Lemma}
\providecommand{\theoremname}{Theorem}

\begin{document}
\title{An Efficient and Globally Optimal Algorithm for Nonconvex QCQP with
One Equality Constraint}
\author{Licheng~Zhao, Rui~Zhou, and Wenqiang~Pu\thanks{Licheng Zhao, Rui Zhou, and Wenqiang Pu are with Shenzhen Research
Institute of Big Data, Shenzhen 518172, China (email: \{zhaolicheng,
wpu, rui.zhou\}@sribd.cn). }}
\maketitle
\begin{abstract}
In this paper, we concentrate on a particular category of quadratically
constrained quadratic programming (QCQP): nonconvex QCQP with one
equality constraint. This type of QCQP problem optimizes a quadratic
objective under a fixed second-order cost and has various engineering
applications. It often serves as a subproblem in an iterative algorithm
framework. However, the development of a high-quality and efficient
solution remains an open problem in the existing literature. Traditionally,
the Semidefinite Relaxation (SDR) technique is applied for an optimal
solution with a prohibitively high order of time complexity. To improve
computational efficiency, we propose a fast and non-iterative algorithm
to reach a globally optimal solution. This algorithm consists of two
consecutive stages: Simultaneous Diagonalization (SD) and Bisection
Search (BS). The SD stage decouples the original problem through an
affine mapping and the BS stage finds the optimal Lagrange multiplier
by solving an equation induced from first- and second-order Karush-Kuhn-Tucker
(KKT) conditions. In addition, we enrich the proposed algorithm with
further extensions on the problem structure, namely, rank-deficient
parameter, indefiniteness, constraint augmentation, and matrix-format
variable. Numerical simulations show that the proposed algorithm achieves
good numerical performance in terms of constraint satisfaction, optimality
gap, and computational time, and scales to problem sizes at least
ten times those supported by the traditional benchmarks.
\end{abstract}


\vspace{-0.4cm}

\section{Introduction}

\IEEEPARstart{Q}{uadratically} constrained quadratic programming
(QCQP) belongs to the family of quadratic optimization problems. In
a QCQP, the objective and constraint functions are both in a quadratic
form. This type pf optimization problem is endowed with ubiquitous
applicability in signal processing, like the max-cut problem \citep{mohar1990eigenvalues},
quadratic knapsack problem \citep{gallo1980quadratic}, multicast
transmit beamforming \citep{sidiropoulos2006transmit,karipidis2007far},
robust receive beamforming \citep{gershman2010convex}, and radar
waveform design \citep{de2008code}. In this paper, we focus on a
special type of QCQP with a single equality constraint: \vspace{-0.3cm}
\begin{equation}
\begin{aligned} & \underset{\mathbf{x}\in\mathbb{R}^{N}}{\mathsf{minimize}} &  & \mathbf{x}^{T}\mathbf{A}_{0}\mathbf{x}+2\mathbf{b}_{0}^{T}\mathbf{x}+c_{0}\\
 & \mathsf{subject\thinspace to} &  & \mathbf{x}^{T}\mathbf{A}_{1}\mathbf{x}+2\mathbf{b}_{1}^{T}\mathbf{x}+c_{1}=0.
\end{aligned}
\label{eq: QCQP eq real}
\end{equation}
This QCQP includes the generalized eigenvalue problem as a special
case \citep{jolliffe2002principal,hotelling1936relations}. The underlying
motivation is to optimize a quadratic target under a fixed second-order
cost, like power budget in wireless communications and volatility
risk in financial engineering. This problem can be interpreted as
a boundary-constrained Trust-Region subproblem \citep{more1993generalizations}
and can serve as a subproblem of an iterative algorithm framework,
like Majorization Minimization (MM) \citep{sun2016majorization} or
Block Coordinate Descent (BCD) \citep[Chap. 2.7]{bertsekas1997nonlinear}. 

\vspace{-0.4cm}

\subsection{Related Works }

A general (real-valued) QCQP can be expressed as: \vspace{-0.2cm}
\begin{equation}
\begin{aligned} & \underset{\mathbf{x}\in\mathbb{R}^{N}}{\mathsf{minimize}} &  & \mathbf{x}^{T}\mathbf{A}_{0}\mathbf{x}+2\mathbf{b}_{0}^{T}\mathbf{x}+c_{0}\\
 & \mathsf{subject\thinspace to} &  & \mathbf{x}^{T}\mathbf{A}_{i}\mathbf{x}+2\mathbf{b}_{i}^{T}\mathbf{x}+c_{i}\leq0,i=1,2,\ldots,I,
\end{aligned}
\label{eq:general QCQP}
\end{equation}
where $N$ is the variable dimension, $I$ is the number of constraint
functions, $\mathbf{A}_{0}$, $\mathbf{A}_{1}$, ..., $\mathbf{A}_{I}$
are $N\times N$ symmetric matrices, $\mathbf{b}_{0}$, $\mathbf{b}_{1}$,
..., $\mathbf{b}_{I}$ are $N\times1$ column vectors, and $c_{0}$,
$c_{1}$, ..., $c_{I}$ are real scalars. For convenient discussion,
we assume the existence of a solution. \citet{park2017general} have
offered a comprehensive overview of QCQP problems. When $\mathbf{A}_{0}$,
$\mathbf{A}_{1}$, ..., $\mathbf{A}_{I}$ are all positive semidefinite,
QCQP \eqref{eq:general QCQP} becomes a convex problem and its global
optimal solution can be attained through numerous efficient algorithms,
like the interior-point method \citep{boyd2004convex}. When $\mathbf{A}_{0}$
or any $\mathbf{A}_{i}$ has negative eigenvalues, the QCQP problem
turns nonconvex and is known to be NP-hard in general \citep{d2003relaxations}.
Even if $\mathbf{A}_{i}$'s are all zero matrices, the degenerated
QCQP (known as Quadratic Programming, short for QP) is NP-complete
\citep{vavasis1990quadratic}. In the nonconvex scenario, even feasible
point pursuit is nontrivial and the authors of \citep{mehanna2014feasible,konar2017first}
put forward several iterative algorithms based on Successive Convex
Approximation (SCA) to tackle this problem. In short, convex QCQPs
are tractable with global optimality guarantee while the optimal solutions
of the nonconvex counterparts can hardly be acquired within polynomial
running time. 

However, the magical technique of Semidefinite Relaxation (SDR) \citep{al1995relaxation,kim2003exact,luo2010semidefinite}
made a breakthrough in solving nonconvex QCQPs. Define $\mathbf{X}=\mathbf{x}\mathbf{x}^{T}$
and QCQP \eqref{eq:general QCQP} can be rewritten as
\begin{equation}
\begin{aligned} & \underset{\mathbf{X}\in\mathbb{R}^{N\times N},\thinspace\mathbf{x}\in\mathbb{R}^{N}}{\mathsf{minimize}} &  & \mathrm{Tr}\left(\mathbf{A}_{0}\mathbf{X}\right)+2\mathbf{b}_{0}^{T}\mathbf{x}+c_{0}\\
 & \mathsf{\hspace{0.4cm}subject\thinspace to} &  & \mathrm{Tr}\left(\mathbf{A}_{i}\mathbf{X}\right)+2\mathbf{b}_{i}^{T}\mathbf{x}+c_{i}\leq0,i=1,2,\ldots,I,\\
 &  &  & \mathbf{X}=\mathbf{x}\mathbf{x}^{T}.
\end{aligned}
\label{eq:QCQP X}
\end{equation}
The SDR technique performs a relaxation on the equality constraint
$\mathbf{X}=\mathbf{x}\mathbf{x}^{T}$, leading to a convex positive-semidefinite-cone
(PSD-cone) constraint $\mathbf{X}-\mathbf{x}\mathbf{x}^{T}\succeq\mathbf{0}$.
This PSD-cone constraint has a Schur complement reformulation and
QCQP \eqref{eq:QCQP X} is relaxed into 
\begin{equation}
\begin{aligned} & \underset{\mathbf{X}\in\mathbb{R}^{N\times N},\thinspace\mathbf{x}\in\mathbb{R}^{N}}{\mathsf{minimize}} &  & \mathrm{Tr}\left(\mathbf{A}_{0}\mathbf{X}\right)+2\mathbf{b}_{0}^{T}\mathbf{x}+c_{0}\\
 & \mathsf{\hspace{0.4cm}subject\thinspace to} &  & \mathrm{Tr}\left(\mathbf{A}_{i}\mathbf{X}\right)+2\mathbf{b}_{i}^{T}\mathbf{x}+c_{i}\leq0,i=1,2,\ldots,I,\\
 &  &  & \left[\begin{array}{cc}
\mathbf{X} & \mathbf{x}\\
\mathbf{x}^{T} & 1
\end{array}\right]\succeq\mathbf{0},
\end{aligned}
\label{eq:QCQP relax}
\end{equation}
which is a Semidefinite Programming (SDP), a typical convex problem.
The computational complexity of solving an SDP is $\mathcal{O}\left(N^{6.5}\right)$
\citep{vandenberghe2005interior,vandenberghe1996semidefinite}. Under
special circumstances, SDR is able to achieve tightness, or equivalently,
the solution $\mathbf{X}$ to QCQP \eqref{eq:QCQP relax} satisfies
the rank-one property \citep{pu2018optimal}. In particular, SDR conditionally
enforces a zero duality gap for nonconvex QCQPs when they are homogeneous
(i.e., $\mathbf{b}_{0}=\mathbf{b}_{i}=\mathbf{0}$, $c_{0}=c_{i}=0$,
$\forall i$) and the number of constraints does not exceed $2$ for
real-valued problems and $3$ for complex-valued problems. Inhomogeneous
QCQP \eqref{eq:general QCQP} can be homogenized via variable augmentation
\citep{luo2010semidefinite}. Further results for tight SDR are elaborated
in \citep{huang2009rank}. In case SDR is not tight, one practical
approach to suboptimal QCQP solutions is randomization \citep{d2003relaxations}.
Besides that, \citet{huang2016consensus} applied Consensus-ADMM for
general QCQPs, and \citet{konar2017fast} used first-order methods
for the max-min nonconvex QCQP. 

Parallel to SDR, Lagrangian relaxation is another alternative, which
is to solve the dual problem of QCQP. The optimal value of the dual
problem provides a lower bound for the original QCQP in general. But
specially, when a QCQP (not necessarily convex) has only one constraint,
as in
\begin{equation}
\begin{aligned} & \underset{\mathbf{x}\in\mathbb{R}^{N}}{\mathsf{minimize}} &  & \mathbf{x}^{T}\mathbf{A}_{0}\mathbf{x}+2\mathbf{b}_{0}^{T}\mathbf{x}+c_{0}\\
 & \mathsf{subject\thinspace to} &  & \mathbf{x}^{T}\mathbf{A}_{1}\mathbf{x}+2\mathbf{b}_{1}^{T}\mathbf{x}+c_{1}\leq0,
\end{aligned}
\label{eq:QCQP ineq}
\end{equation}
Lagrangian relaxation is tight and a zero duality gap is guaranteed.
Lagrangian relaxation tightness is justified by the S-procedure theorem
\citep[Appendix B.2]{boyd2004convex}. Early researches on single-constraint
QCQPs include \citep{gander1978linear,golub1991quadratically,ben1996hidden}.
This special type of QCQP problem readily finds its application in
the Trust-Region subproblem of an iterative algorithm framework \citep{more1993generalizations,conn2000trust,yuan2000review,yuan2015recent}.
Although single-constraint QCQPs enjoy equivalence to their dual problems
in terms of optimal value, regardless of convexity, the solving algorithm
is still iterative. The dual problem of a QCQP is expressed as an
SDP \citep[Appendix B.1]{boyd2004convex}, so the overall computational
complexity is still $\mathcal{O}\left(N^{6.5}\right)$. This order
of complexity is prohibitively high and not affordable for large-scale
problems. 

Recent advances in single-constraint QCQP problems are elaborated in
\citep{feng2012duality,adachi2019eigenvalue,taati2020local,song2023local}.
Notably, \citet{song2023local} proved the necessity of the standard
second-order optimality condition for a local non-global minimizer.
\citet{adachi2019eigenvalue} tackled the single-constraint QCQP with
an efficient algorithm, of complexity $\mathcal{O}\left(N^{3}\right)$.
This algorithm was developed from a set of necessary and sufficient
conditions by \citet{more1993generalizations} on top of the Slater\textquoteright s
condition. Those conditions are actually first- and second-order Karush-Kuhn-Tucker
(KKT) conditions \citep{boyd2004convex,nocedal1999numerical}. Primal-dual
methods well suit the single-inequality-constraint QCQP problem because
the underlying strong duality only requires the Slater\textquoteright s
condition. The algorithm developed in \citep{adachi2019eigenvalue}
involves Eigenvalue Decomposition (EVD), Singular Value Decomposition
(SVD), and linear system solving. Note that this algorithm is not
iterative, so it should be regarded as an analytic solution. However,
the bottleneck of \citeauthor{adachi2019eigenvalue}'s algorithm is
to solve a $\left(2N+1\right)$-dimensional generalized eigenvalue
problem for an extremal eigenpair. This submodule is implemented with
a software package named ``ARPACK'' \citep{lehoucq1998arpack} (this
package is embedded in the Matlab function ``eigs'') and may fail
to converge in face of close-to-singularity situations. 

QCQP \eqref{eq: QCQP eq real} confines the feasible set on the border,
which is a restricted variant of problem \eqref{eq:QCQP ineq}. This
problem has been extensively studied in the existing literature \citep{moler1973algorithm,zhao2018mean,yu2016alternating,mao2024joint}.
In particular, when $\mathbf{b}_{0}=\mathbf{b}_{1}=\mathbf{0}$, $c_{0}=0$,
and $c_{1}=-1$, this special type of QCQP is reduced to the famous
generalized eigenvalue problem. In financial engineering or wireless
communications, we may want to maximize a utility function subject
to a fixed volatility risk or constant power budget constraint. A
fixed volatility risk ensures reasonable modeling of portfolio optimization.
A constant power budget helps avoid energy waste and build green and
low-carbon communication systems. If the utility function is successively
approximated with a quadratic surrogate function, then we need to
solve a series of single-equality-constraint QCQPs. For example, in
mean-reverting portfolio design \citep{zhao2018mean}, the authors
applied the MM framework and the majorized subproblem is a single-equality-constraint
QCQP. This QCQP is inhomogeneous in the case of a net-budget portfolio
and homogeneous for a zero-budget portfolio. Another example is the
millimeter wave multiple-input--multiple-output (MIMO) system design
in wireless communications \citep{yu2016alternating}. The hybrid
precoder optimization problem is addressed in an alternating minimization
(i.e., BCD) manner and one alternating subproblem is a single-equality-constraint
QCQP with a matrix-format variable, cf. \citep[eq. (34)]{yu2016alternating}.
\citeauthor{yu2016alternating} proposed to solve this QCQP via matrix
vectorization and then SDR, so the high order complexity could be
an algorithmic bottleneck. Single-equality-constraint QCQPs also show
up in the application of angles of arrival (AOA) and differential
time delay (DTD) localization. Interested readers may refer to \citep{mao2024joint}
for further information. 

So far, literature on efficient algorithms for solving \eqref{eq: QCQP eq real}
is still scarce. In this problem, primal-dual methods are not readily
applicable due to weak duality. The first- and second-order KKT conditions
as in \citep{more1993generalizations} are no longer sufficient but
merely necessary. From the numerical perspective, a simple gradient
projection method cannot be easily implemented because of the nontrivial
projection operation onto an equality-constrained quadratic set. Hence,
the global solution of \eqref{eq: QCQP eq real} still relies on the
traditional SDR technique. In this paper, we intend to solve several
classes of single-equality-constraint QCQPs with better efficiency
than SDR. Representative QCQP problem types and solution methods are
summarized in Table \ref{tab:ref summary} for easy reference. 

\begin{table*}[tbh]
\centering
\begin{threeparttable}
\caption{A summary of representative works in QCQP optimization.}
\label{tab:ref summary}
\begin{centering}
\begin{tabular}{|c|c|c|}
\hline 
QCQP Problem Types & Solution Method & Complexity\tabularnewline
\hline 
\multirow{2}{*}{General QCQP } & Convex: interior-point method \citep{boyd2004convex} & $\mathcal{O}\left(N^{3}\right)$\tabularnewline
\cline{2-3} \cline{3-3} 
 & Nonconvex: SDR and randomization \citep{park2017general,luo2010semidefinite}  & $\mathcal{O}\left(N^{6.5}\right)$\tabularnewline
\hline 
\multirow{2}{*}{Single-inequality-constraint QCQP } & Lagrangian relaxation or S-procedure \citep{boyd2004convex} and SDR
\citep{luo2010semidefinite} & $\mathcal{O}\left(N^{6.5}\right)$ \tabularnewline
\cline{2-3} \cline{3-3} 
 & An EVD-based method \citep{adachi2019eigenvalue} & $\mathcal{O}\left(N^{3}\right)$\tabularnewline
\hline 
Single-equality-constraint QCQP  & SDR \citep{luo2010semidefinite} & $\mathcal{O}\left(N^{6.5}\right)$\tabularnewline
\hline 
\end{tabular}
\par\end{centering}
\centering{}
\end{threeparttable}
\end{table*}

\subsection{Contribution}

The major contributions of this paper are listed as follows:
\begin{enumerate}
\item We propose a non-iterative globally optimal algorithm for the nonconvex
single-equality-constraint QCQP problems. This algorithm is computationally
efficient and designed with a two-stage scheme: Simultaneous Diagonalization
(SD) and Bisection Search (BS). The SD stage decouples the original
problem through an affine mapping. The BS stage finds the optimal
Lagrange multiplier of the first- and second-order KKT system. 
\item Besides the standard case, we investigate into a few extensions for
further exploration, including rank deficiency, indefiniteness, linear
constraint augmentation, and matrix-format optimization variable.
To clearly demonstrate the QCQPs within our solvability scope, we
present the problem characteristics in the table below. 
\begin{table}[tbh]
\centering{}\scalebox{1.07}{%
\begin{tabular}{|c|c|}
\hline 
 & Problem Characteristics\tabularnewline
\hline 
Standard case & $\mathbf{A}_{0}$ symmetric and $\mathbf{A}_{1}\succ\mathbf{0}$\tabularnewline
\hline 
Extension I: & $\mathbf{A}_{0}\succeq\mathbf{0}$,\tabularnewline
Rank Deficiency & $\mathbf{A}_{1}\succeq\mathbf{0}$ and $\mathbf{A}_{1}\nsucc\mathbf{0}$\tabularnewline
\hline 
Extension II: & $\mathbf{A}_{0}\succ\mathbf{0}$,\tabularnewline
Indefiniteness & $\mathbf{A}_{1}$ symmetric, indefinite, and full-rank\tabularnewline
\hline 
Extension III: & Additional linear equalities\tabularnewline
Constraint Augmentation & $\mathbf{A}_{0}$ symmetric and $\mathbf{A}_{1}\succ\mathbf{0}$\tabularnewline
\hline 
Extension IV: & From $\mathbf{x}$ to $\mathbf{X}$\tabularnewline
Matrix-format Variable & $\mathbf{A}_{0}$ symmetric and $\mathbf{A}_{1}\succ\mathbf{0}$\tabularnewline
\hline 
\end{tabular}}
\end{table}
\item The proposed algorithm achieves good numerical performance in the
following four aspects: constraint satisfaction, optimality gap, computational
time, and scalability. In general, the constraint satisfaction levels
are universally below $10^{-5}$ and the optimality gaps are lower
than $10^{-4}$. The computational time is at least two orders of
magnitude shorter than the compared benchmarks. The computation time
is reduced by at least two orders of magnitude compared to the existing
benchmarks. With regard to scalability, the proposed algorithm can
handle problem sizes at least ten times as large as those affordable
by an off-the-shelf SDP solver. Specifically, in the scenario of complex-valued
matrix-format QCQP, the proposed algorithm is able to produce a lower
objective value as well as a higher computational speed while maintaining
feasibility. 
\end{enumerate}

\subsection{Organization and Notation }

The rest of the paper is organized as follows. In Section \ref{sec:Globally-Optimal-Solution},
we develop the globally optimal solution to the standard case of nonconvex
QCQP with one equality constraint. In Sections \ref{sec:Quadratic-Matrix-Extension},
\ref{sec: Indefiniteness}, \ref{sec:Constraint-Augmentation-Extension},
and \ref{sec:Matrix-Format-Extension}, we explore several extensions
based on the standard case, namely, rank-deficient parameter, indefiniteness,
constraint augmentation, and matrix-format variable. Finally, Section
\ref{sec:Numerical-Simulations} displays numerical simulations, and
the conclusions are drawn in Section \ref{sec:Conclusions}. 

The following notation is adopted. Boldface upper-case letters represent
matrices, boldface lower-case letters denote column vectors, and standard
lower-case or upper-case letters stand for scalars. $\mathbb{R}$
($\mathbb{C}$) denotes the real (complex) field. $\left|\cdot\right|$
denotes the absolute value for the real case and modulus for the complex
case. $\mathrm{Re}\left[\cdot\right]$ denotes the real part of a
complex number. $\arg\left(\cdot\right)$ stands for the phase of
a complex number. $\mathrm{card}\left(\cdot\right)$ represents the
cardinality of a set. $\left\Vert \cdot\right\Vert _{p}$ denotes
the $\ell_{p}$-norm of a vector. $\nabla\left(\cdot\right)$ represents
the gradient of a multivariate function (the way to derive the complex-valued
gradient follows Euclidean gradient in \citep{absil2009optimization}).
$\odot$ stands for the Hadamard product. $\mathbf{1}$ stands for
the all-one vector, and $\mathbf{I}$ stands for the identity matrix.
$\mathbf{X}^{T}$, $\mathbf{X}^{*}$, $\mathbf{X}^{H}$, $\mathbf{X}^{-1}$,
$\mathbf{X}^{\dagger}$, $\mathcal{R}\left(\mathbf{X}\right)$, $\textrm{Tr}\left(\mathbf{X}\right)$,
and $\textrm{vec}\left(\mathbf{X}\right)$ denote the transpose, complex
conjugate, conjugate transpose, inverse, pseudo-inverse, column space,
trace, and stacking vectorization of $\mathbf{X}$, respectively.
$\mathbf{X}\succeq\mathbf{0}\left(\succ\mathbf{0}\right)$ means that
$\mathbf{X}$ is positive semidefinite (positive definite). $\left\Vert \mathbf{X}\right\Vert _{F}$
is the Frobenius norm of $\mathbf{X}$. $\textrm{Diag}\left(\mathbf{x}\right)$
is a diagonal matrix with $\mathbf{x}$ filling its principal diagonal.
The superscript $\star$ represents the optimal solution to an optimization
problem. Whenever arithmetic operators ($\sqrt{\cdot}$, $\cdot/\cdot$,
$\cdot^{2}$, $\left|\cdot\right|$, $\max$, etc.) are applied to
vectors or matrices, they are elementwise operations. 

\section{Globally Optimal Solution Pursuit \label{sec:Globally-Optimal-Solution}}

Let's recall the single-equality-constraint QCQP \eqref{eq: QCQP eq real}.
We start with the standard case where $\mathbf{A}_{0}$ is symmetric
and $\mathbf{A}_{1}\succ\mathbf{0}$. Since $\mathbf{A}_{1}\succ\mathbf{0}$,
the variable is bounded in every dimension. If the constraint set
is nonempty, an optimal solution must exist. Thus, we further assume
$-\mathbf{b}_{1}^{T}\mathbf{A}_{1}^{-1}\mathbf{b}_{1}+c_{1}\leq0$
to ensure solution existence. In the standard case, the globally optimal
solution scheme consists of two stages: Simultaneous Diagonalization
(SD) and Bisection Search (BS). 

\subsection{Simultaneous Diagonalization\label{subsec:Simultaneous-Diagonalization}}

The first stage of the solution scheme is to perform simultaneous
diagonalization on $\mathbf{A}_{0}$ and $\mathbf{A}_{1}$. Simultaneous
diagonalization finds an affine mapping for the original variable
$\mathbf{x}$ so that it is elementwisely decoupled after transformation.
Simultaneous diagonalization is realized through the following lemma. 
\begin{lem}[\citep{horn2012matrix}]
\label{lem:SD-full}For any symmetric $\mathbf{A}_{0}$ and $\mathbf{A}_{1}\succ\mathbf{0}$,
there exists an invertible matrix $\mathbf{T}$ such that $\mathbf{A}_{0}=\mathbf{T}\boldsymbol{\Lambda}_{0}\mathbf{T}^{T}$
and $\mathbf{A}_{1}=\mathbf{T}\mathbf{T}^{T}$ where $\boldsymbol{\Lambda}_{0}$
is diagonal. 
\end{lem}
\begin{IEEEproof}
Since $\mathbf{A}_{1}\succ\mathbf{0}$, there exists an invertible
square matrix $\mathbf{S}$ such that $\mathbf{S}\mathbf{S}^{T}=\mathbf{A}_{1}$.
Then we perform EVD on $\mathbf{S}^{-1}\mathbf{A}_{0}\mathbf{S}^{-T}$:
\begin{equation}
\mathbf{S}^{-1}\mathbf{A}_{0}\mathbf{S}^{-T}=\mathbf{U}\boldsymbol{\Lambda}_{0}\mathbf{U}^{T}\label{eq:simul diag A0}
\end{equation}
where $\mathbf{U}$ is an unitary matrix and $\boldsymbol{\Lambda}_{0}=\mathrm{Diag}\left(\mathbf{a}\right)$
is diagonal. To this end, we are able to diagonalize $\mathbf{A}_{0}$
and $\mathbf{A}_{1}$ simultaneously: 
\begin{equation}
\begin{cases}
\mathbf{A}_{0}=\mathbf{S}\mathbf{U}\boldsymbol{\Lambda}_{0}\mathbf{U}^{T}\mathbf{S}^{T}=\left(\mathbf{S}\mathbf{U}\right)\boldsymbol{\Lambda}_{0}\left(\mathbf{S}\mathbf{U}\right)^{T}=\mathbf{T}\boldsymbol{\Lambda}_{0}\mathbf{T}^{T}\\
\mathbf{A}_{1}=\mathbf{S}\mathbf{S}^{T}=\mathbf{S}\mathbf{U}\cdot\mathbf{I}\cdot\mathbf{U}^{T}\mathbf{S}^{T}=\left(\mathbf{S}\mathbf{U}\right)\mathbf{I}\left(\mathbf{S}\mathbf{U}\right)^{T}=\mathbf{T}\mathbf{T}^{T}
\end{cases}\label{eq:sim diag}
\end{equation}
with $\mathbf{T}\triangleq\mathbf{S}\mathbf{U}$. Since $\mathbf{S}$
and $\mathbf{U}$ are both full-rank, $\mathbf{T}$ is invertible
as well. 
\end{IEEEproof}
Next, we perform the affine mapping on the original variable $\mathbf{x}$.
The decoupling process is summarized in the theorem below. 
\begin{thm}
\label{thm:eqv recast}The constraint and objective functions in QCQP
\eqref{eq: QCQP eq real} can be equivalently recast as 
\begin{equation}
\begin{cases}
\mathbf{x}^{T}\mathbf{A}_{1}\mathbf{x}+2\mathbf{b}_{1}^{T}\mathbf{x}+c_{1}=\mathbf{y}^{T}\mathbf{y}-c\\
\mathbf{x}^{T}\mathbf{A}_{0}\mathbf{x}+2\mathbf{b}_{0}^{T}\mathbf{x}+c_{0}=\mathbf{y}^{T}\boldsymbol{\Lambda}_{0}\mathbf{y}-2\mathbf{b}^{T}\mathbf{y}+\mathrm{const}
\end{cases}
\end{equation}
where $\mathbf{y}=\mathbf{T}^{T}\mathbf{x}+\mathbf{T}^{-1}\mathbf{b}_{1}$,
$c=\left(\mathbf{T}^{-1}\mathbf{b}_{1}\right)^{T}\left(\mathbf{T}^{-1}\mathbf{b}_{1}\right)-c_{1}\geq0$,
$\mathbf{b}=\boldsymbol{\Lambda}_{0}\mathbf{T}^{-1}\mathbf{b}_{1}-\mathbf{T}^{-1}\mathbf{b}_{0}$,
and $\mathrm{const}$ is an optimization irrelevant constant. 
\end{thm}
\begin{IEEEproof}
See Appendix \ref{sec:Proof-of-Theorem eqv recast} for the detailed
proof. 
\end{IEEEproof}
Under the assumption $-\mathbf{b}_{1}^{T}\mathbf{A}_{1}^{-1}\mathbf{b}_{1}+c_{1}\leq0$,
$c\geq0$ holds true. Hence, the original QCQP \eqref{eq: QCQP eq real}
is eventually reformulated as
\begin{equation}
\begin{aligned} & \underset{\mathbf{y}\in\mathbb{R}^{N}}{\mathsf{minimize}} &  & \mathbf{y}^{T}\boldsymbol{\Lambda}_{0}\mathbf{y}-2\mathbf{b}^{T}\mathbf{y}\\
 & \mathsf{subject\thinspace to} &  & \mathbf{y}^{T}\mathbf{y}-c=0.
\end{aligned}
\end{equation}
Note that $\mathbf{x}$ can be recovered from the inverse mapping
of $\mathbf{y}$. We rewrite the problem by replacing $\boldsymbol{\Lambda}_{0}$
with $\mathrm{Diag}\left(\mathbf{a}\right)$ so that the optimization
variables are completely decoupled: 

\begin{equation}
\begin{aligned} & \underset{\left\{ y_{i}\right\} _{i=1}^{n}}{\mathsf{minimize}} &  & \sum_{i=1}^{N}a_{i}y_{i}^{2}-2b_{i}y_{i}\\
 & \mathsf{subject\thinspace to} &  & \sum_{i=1}^{N}y_{i}^{2}-c=0.
\end{aligned}
\label{eq:scalar QCQP}
\end{equation}

\subsection{Bisection Search \label{subsec:Bisection-Search}}

When $c=0$, the only feasible solution is $y_{i}=0$, $\forall i$,
and we obtain the trivial solution $\mathbf{x}^{\star}=-\mathbf{T}^{-T}\mathbf{T}^{-1}\mathbf{b}_{1}=-\mathbf{A}_{1}^{-1}\mathbf{b}_{1}$.
When $c>0$, the feasible set is nonempty and we can derive the optimal
solution from KKT conditions. Because problem \eqref{eq:scalar QCQP}
is nonconvex, a KKT point is not necessarily optimal but the optimal
solution must lie among the KKT points \citep[Chap. 5]{boyd2004convex}.
We start with the first-order KKT conditions: 1) primal feasibility:
\begin{equation}
\sum_{i=1}^{N}y_{i}^{2}-c=0\label{eq:cond feasible}
\end{equation}
and 2) zero-gradient Lagrangian ($\lambda$ is the Lagrange multiplier):
\begin{equation}
\begin{aligned} & \frac{\partial\mathcal{L}}{\partial y_{i}}=0,\forall i\hspace{0.3cm}\Longrightarrow\\
 & \frac{\partial}{\partial y_{i}}\left(\left(\sum_{i=1}^{N}a_{i}y_{i}^{2}-2b_{i}y_{i}\right)+\lambda\left(\left(\sum_{i=1}^{N}y_{i}^{2}\right)-c\right)\right)=0,\forall i.
\end{aligned}
\label{eq:cond grad}
\end{equation}
Condition \eqref{eq:cond grad} gives rise to the optimal $y_{i}$'s:
\begin{equation}
y_{i}=b_{i}/\left(a_{i}+\lambda\right).
\end{equation}
The primal feasibility condition is reduced to 
\begin{equation}
\sum_{i=1}^{N}b_{i}^{2}/\left(a_{i}+\lambda\right)^{2}=c.
\end{equation}
It is straightforward that $\lambda\neq-a_{i}$, $\forall i$. 

\begin{figure}[tbh]
\begin{centering}
\includegraphics[width=0.65\columnwidth]{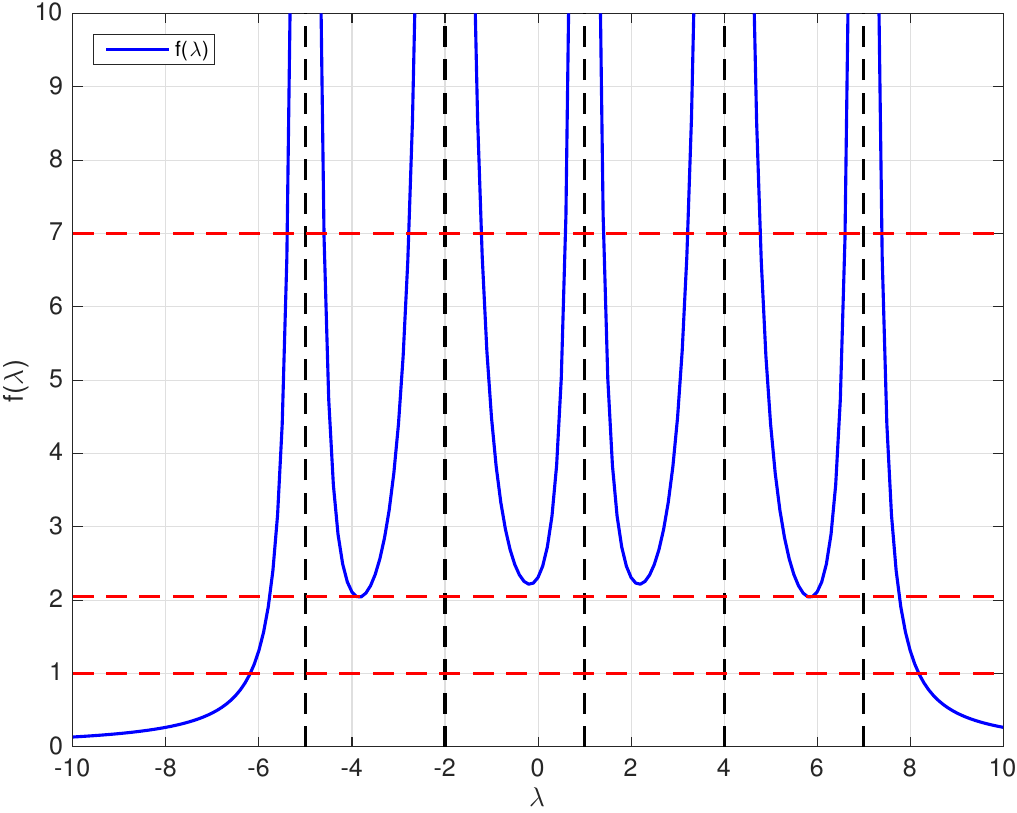}
\par\end{centering}
\caption{\label{fig:An-illustration-plot}An illustration plot of $f\left(\lambda\right)=\frac{1}{\left(\lambda+5\right)^{2}}+\frac{4}{\left(\lambda+2\right)^{2}}+\frac{1}{\left(\lambda-1\right)^{2}}+\frac{4}{\left(\lambda-4\right)^{2}}+\frac{1}{\left(\lambda-7\right)^{2}}$.}
\end{figure}

We denote the function on the left hand side as $f\left(\lambda\right)$
\begin{equation}
f\left(\lambda\right)=\sum_{i=1}^{N}b_{i}^{2}/\left(a_{i}+\lambda\right)^{2}.\label{eq:f lambda}
\end{equation}
Before solving the equation $f\left(\lambda\right)=c$, let's examine
the shape of $f\left(\lambda\right)$. Below is an illustration plot
of $f\left(\lambda\right)$ where $N=5$, $a_{1}=5$, $a_{2}=2$,
$a_{3}=-1$, $a_{4}=-4$, $a_{5}=-7$, $b_{1}=b_{3}=b_{5}=1$, and
$b_{2}=b_{4}=2$. As can be seen in Figure \ref{fig:An-illustration-plot},
$f\left(\lambda\right)$ is continuous almost everywhere except for
a few essential discontinuity points and different choices of $c$
lead to different number of roots in $f\left(\lambda\right)=c$. 
 The monotonicity behavior of $f\left(\lambda\right)$ indicates multiple
KKT points in the first-order equation system. To figure out the optimal
solution from the first-order KKT candidates, we need to introduce
the second-order KKT condition \citep{nocedal1999numerical}: 
\begin{equation}
\begin{aligned} & \frac{\partial^{2}\mathcal{L}}{\partial y_{i}^{2}}\geq0,\forall i\hspace{0.3cm}\\
\Longrightarrow\hspace{0.2cm} & 2\left(a_{i}+\lambda\right)\geq0,\forall i\\
\Longrightarrow\hspace{0.2cm} & \lambda\geq-\min\left(\left\{ a_{i}\right\} \right)
\end{aligned}
\label{eq:cond hess}
\end{equation}
This is also a necessary condition for optimality. The intuition
is that any global minimum of the optimization problem must be a local
minimum on the Lagrangian function. Because $\lambda=-\min\left(a_{i}\right)$
is a singular point of $f\left(\lambda\right)$, the optimal $\lambda$
belongs to $\left(-\min\left(\left\{ a_{i}\right\} \right),+\infty\right)$. 

After introducing the second-order KKT condition, we continue with
the equation $f\left(\lambda\right)=c$. Within $\lambda\in\left(-\min\left(\left\{ a_{i}\right\} \right),+\infty\right)$,
the first derivative of $f\left(\lambda\right)$ is non-positive: 

\begin{equation}
f^{\prime}\left(\lambda\right)=-\sum_{i=1}^{N}2b_{i}^{2}\left(a_{i}+\lambda\right)^{-3}\leq0.
\end{equation}
Therefore, $f\left(\lambda\right)$ is monotonically decreasing in
the specified interval. In additional, when $\lambda\rightarrow\left[-\min\left(\left\{ a_{i}\right\} \right)\right]_{+}$,
$f\left(\lambda\right)\rightarrow+\infty$; and when $\lambda\rightarrow+\infty$,
$f\left(\lambda\right)\rightarrow0$. Therefore, for any $c>0$, there
exist a unique root to $f\left(\lambda\right)=c$ and due to monotonicity,
this root can be found with bisection search, denoted by $\lambda^{\star}$.
The optimal solution $\mathbf{y}^{\star}$ is $\mathbf{y}^{\star}=\mathbf{b}/\left(\mathbf{a}+\lambda^{\star}\mathbf{1}\right)$
and the optimal $\mathbf{x}^{\star}$ is 
\begin{equation}
\mathbf{x}^{\star}=\mathbf{T}^{-T}\mathbf{y}^{\star}-\mathbf{T}^{-T}\mathbf{T}^{-1}\mathbf{b}_{1}=\mathbf{T}^{-T}\mathbf{y}^{\star}-\mathbf{A}_{1}^{-1}\mathbf{b}_{1}.
\end{equation}

The second stage of the solution scheme is to solve the first- and
second-order KKT system. This system turns out to have a unique KKT
point, so it automatically generates the globally optimal solution.
The entire process is shortly named as ``SD-BS'' and the complete
procedure is summarized in Algorithm \ref{alg:Globally-optimal-solution}.

\begin{algorithm}[tbh]
\caption{\label{alg:Globally-optimal-solution}Globally optimal solution to
single-equality-constraint QCQPs: Simultaneous Diagonalization - Bisection
Search (SD-BS). }

\begin{algorithmic}[1] 

\REQUIRE QCQP parameters: $\mathbf{A}_{0}$, $\mathbf{b}_{0}$, $\mathbf{A}_{1}$,
$\mathbf{b}_{1}$, and $c_{1}$. 

\STATE Perform EVD on $\mathbf{A}_{1}$: $\mathbf{A}_{1}=\mathbf{U}_{1}\boldsymbol{\Lambda}_{1}\mathbf{U}_{1}^{T}$; 

\STATE Compute $\mathbf{S}^{-1}$: $\mathbf{S}^{-1}=\boldsymbol{\Lambda}_{1}^{-1/2}\mathbf{U}_{1}^{T}$;

\STATE Perform EVD on $\mathbf{S}^{-1}\mathbf{A}_{0}\mathbf{S}^{-T}$:
$\mathbf{S}^{-1}\mathbf{A}_{0}\mathbf{S}^{-T}=\mathbf{U}\boldsymbol{\Lambda}_{0}\mathbf{U}^{T}$
with $\boldsymbol{\Lambda}_{0}=\mathrm{Diag}\left(\mathbf{a}\right)$;

\STATE Compute $\mathbf{T}^{-1}$: $\mathbf{T}^{-1}=\mathbf{U}^{T}\mathbf{S}^{-1}$;

\STATE Compute $\mathbf{b}$ and $c$: $\mathbf{b}=\boldsymbol{\Lambda}_{0}\left(\mathbf{T}^{-1}\mathbf{b}_{1}\right)-\mathbf{T}^{-1}\mathbf{b}_{0}$
and $c=\left(\mathbf{T}^{-1}\mathbf{b}_{1}\right)^{T}\left(\mathbf{T}^{-1}\mathbf{b}_{1}\right)-c_{1}$;

\STATE In $\left(-\min\left(\left\{ a_{i}\right\} \right),+\infty\right)$,
find $\lambda^{\star}$ such that $f\left(\lambda^{\star}\right)=c$
through bisection search; 

\STATE Compute $\mathbf{y}^{\star}$: $\mathbf{y}^{\star}=\mathbf{b}/\left(\mathbf{a}+\lambda^{\star}\mathbf{1}\right)$;

\STATE Compute $\mathbf{x}^{\star}$: $\mathbf{x}^{\star}=\mathbf{T}^{-T}\left(\mathbf{y}^{\star}-\mathbf{T}^{-1}\mathbf{b}_{1}\right)$. 

\end{algorithmic} 
\end{algorithm}

\subsection{Computational Complexity}

The dominant computational cost of Algorithm SD-BS results from the
EVD and matrix-matrix multiplication operation, both of complexity
$\mathcal{O}\left(N^{3}\right)$. Algorithm SD-BS involves two EVD's
and two matrix-matrix multiplications if performed in the correct
calculation order. Therefore, the worst-case computational complexity
is $\mathcal{O}\left(N^{3}\right)$. 

\section{Extension I: Rank Deficiency\label{sec:Quadratic-Matrix-Extension}}

From this section onward, we would like to extend the standard case.
In the standard case, all variable components are bounded. But what
if the components are only partially bounded? This question leads
to the first extension on the rank deficiency issue. Precisely, we
assume $\mathbf{A}_{0}\succeq\mathbf{0}$ and $\mathbf{A}_{1}\succeq\mathbf{0}$
but $\mathbf{A}_{1}\nsucc\mathbf{0}$. We will show that the ``SD-BS''
scheme is still applicable in the extension scenario. First we carry
out simultaneous diagonalization on $\mathbf{A}_{0}$ and $\mathbf{A}_{1}$,
as in Lemma \ref{lem:SD-def}. 
\begin{lem}[\citep{bernstein2009matrix}]
\label{lem:SD-def}For any $\mathbf{A}_{0}\succeq\mathbf{0}$, $\mathbf{A}_{1}\succeq\mathbf{0}$,
and $\mathrm{rank}\left(\mathbf{A}_{1}\right)=r<N$, there exists
an invertible matrix $\mathbf{T}$ such that $\mathbf{A}_{0}=\mathbf{T}\boldsymbol{\Lambda}_{0}\mathbf{T}^{T}$
and $\mathbf{A}_{1}=\mathbf{T}\left[\begin{array}{cc}
\mathbf{I} & \mathbf{0}\\
\mathbf{0} & \mathbf{0}
\end{array}\right]\mathbf{T}^{T}$ where $\boldsymbol{\Lambda}_{0}$ is diagonal and the size of $\mathbf{I}$
is $r\times r$. 
\end{lem}
\begin{IEEEproof}
See Appendix \ref{sec:Proof-of-Lemma SD-def} for the detailed proof. 
\end{IEEEproof}
The affine mapping from $\mathbf{x}$ to $\mathbf{y}$ follows Theorem
\ref{thm:eqv recast} in Section \ref{subsec:Simultaneous-Diagonalization}.
We denote $\mathbf{b}_{\mathrm{c}}$ as $\left[\begin{array}{cc}
\mathbf{0} & \mathbf{0}\\
\mathbf{0} & \mathbf{I}
\end{array}\right]\mathbf{T}^{-1}\mathbf{b}_{1}$ and redefine $c=\left(\mathbf{T}^{-1}\mathbf{b}_{1}\right)^{T}\left(\mathbf{T}^{-1}\mathbf{b}_{1}\right)+\left(\mathbf{T}^{-1}\mathbf{b}_{1}\right)^{T}\left[\begin{array}{cc}
\mathbf{0} & \mathbf{0}\\
\mathbf{0} & \mathbf{I}
\end{array}\right]\left(\mathbf{T}^{-1}\mathbf{b}_{1}\right)-c_{1}$. The decoupled QCQP formulation is recast as (the first $r$ elements
of $\mathbf{b}_{\mathrm{c}}$ are zero)
\begin{equation}
\begin{aligned} & \underset{\left\{ y_{i}\right\} _{i=1}^{n}}{\mathsf{minimize}} &  & \sum_{i=1}^{N}a_{i}y_{i}^{2}-2b_{i}y_{i}\\
 & \mathsf{subject\thinspace to} &  & \sum_{i=1}^{r}y_{i}^{2}+2\sum_{i=r+1}^{N}b_{\mathrm{c},i}y_{i}-c=0,
\end{aligned}
\label{eq:scalar QCQP r}
\end{equation}
where $a_{i}\geq0$, $\forall i$. 

Next we move on to the stage of bisection search. We collect the indices
of zero-valued $a_{i}$'s, where $r+1\leq i\leq N$, into a set $\mathcal{I}_{z}=\left\{ i|a_{i}=0,r+1\leq i\leq N\right\} $.
The solution to problem \eqref{eq:scalar QCQP r} exists if and only
if
\begin{equation}
\lambda b_{\mathrm{c},i}-b_{i}=0,\forall i\in\mathcal{I}_{z}.\label{eq:over-determined linear system}
\end{equation}
This condition is interpreted as null space elimination in \citep[Chap. 5.3.1]{adachi2019eigenvalue}.
By exploring the first-order KKT conditions, we obtain the optimal
$y_{i}$'s as 

\begin{equation}
y_{i}=\begin{cases}
b_{i}/\left(a_{i}+\lambda\right) & 1\leq i\leq r\\
\left(b_{i}-\lambda b_{\mathrm{c},i}\right)/a_{i} & r+1\leq i\leq N,i\notin\mathcal{I}_{z}\\
\mathrm{arbitrary} & i\in\mathcal{I}_{z}.
\end{cases}
\end{equation}
When $\mathcal{I}_{z}\neq\phi$, the optimal $\lambda$ is simply
solved from \eqref{eq:over-determined linear system}. When $\mathcal{I}_{z}=\phi$,
we refer to the equality condition of \eqref{eq:scalar QCQP r} and
solve for $\lambda$ from $f\left(\lambda\right)=c$ where 
\begin{equation}
f\left(\lambda\right)=\sum_{i=1}^{r}b_{i}^{2}/\left(a_{i}+\lambda\right)^{2}+2\sum_{i=r+1}^{N}\left(b_{\mathrm{c},i}b_{i}-\lambda b_{\mathrm{c},i}^{2}\right)/a_{i}.
\end{equation}
The newly defined $f\left(\lambda\right)$ may also have multiple
roots on $f\left(\lambda\right)=c$, so the introduction of second-order
KKT condition is still necessary:
\begin{equation}
\lambda\geq-\min_{i=1,\ldots,r}\left(\left\{ a_{i}\right\} \right).
\end{equation}
In $\left(-\min_{i=1,\ldots,r}\left(\left\{ a_{i}\right\} \right),+\infty\right)$,
$f^{\prime}\left(\lambda\right)\leq0$ still holds true:
\begin{equation}
f^{\prime}\left(\lambda\right)=-\sum_{i=1}^{r}2b_{i}^{2}\left(a_{i}+\lambda\right)^{-3}-2\sum_{i=r+1}^{N}\left(b_{\mathrm{c},i}^{2}/a_{i}\right)\leq0,
\end{equation}
where $a_{i}>0$ on $r+1\leq i\leq N$ and $\mathcal{I}_{z}=\phi$.
Additionally, $f\left(\lambda\right)\in\mathbb{R}$ in the specified
interval. Therefore, the optimal $\lambda$ can be attained from bisection
search. 

To sum up, in the rank deficient scenario, we solve QCQP \eqref{eq: QCQP eq real}
under the assumptions $\mathbf{A}_{0}\succeq\mathbf{0}$ and $\mathbf{A}_{1}\succeq\mathbf{0}$
but $\mathbf{A}_{1}\nsucc\mathbf{0}$. To ensure solution existence,
condition \eqref{eq:over-determined linear system} must be satisfied.
The solution procedure is very analogous to Algorithm \ref{alg:Globally-optimal-solution}
up to minor modifications.

\section{Extension II: Indefiniteness\label{sec: Indefiniteness}}

After considering the rank deficiency issue, we would like to take
one step further in this section and make the feasible set unbounded.
All variable components are allowed to go to infinity. This extension
can be realized by the indefinite property of $\mathbf{A}_{1}$. To
ensure solution existence, we additionally require the objective to
be strictly convex. Thus, the assumptions of this extension scenario
are summarized as $\mathbf{A}_{0}\succ\mathbf{0}$, and $\mathbf{A}_{1}$
is symmetric, indefinite, and full-rank. 

We start with the simultaneous diagonalization stage. By a slight
manipulation of Lemma \ref{lem:SD-full}, we have
\begin{equation}
\begin{cases}
\mathbf{A}_{0}=\mathbf{S}\mathbf{S}^{T}=\mathbf{S}\mathbf{U}\cdot\mathbf{I}\cdot\mathbf{U}^{T}\mathbf{S}^{T}=\left(\mathbf{S}\mathbf{U}\right)\mathbf{I}\left(\mathbf{S}\mathbf{U}\right)^{T}=\mathbf{T}\mathbf{T}^{T}\\
\mathbf{A}_{1}=\mathbf{S}\mathbf{U}\boldsymbol{\Lambda}_{1}\mathbf{U}^{T}\mathbf{S}^{T}=\left(\mathbf{S}\mathbf{U}\right)\boldsymbol{\Lambda}_{1}\left(\mathbf{S}\mathbf{U}\right)^{T}=\mathbf{T}\boldsymbol{\Lambda}_{1}\mathbf{T}^{T},
\end{cases}
\end{equation}
where the EVD of $\mathbf{S}^{-1}\mathbf{A}_{1}\mathbf{S}^{-T}$ is
given by $\mathbf{U}\boldsymbol{\Lambda}_{1}\mathbf{U}^{T}$ with
$\mathbf{U}$ being unitary and $\boldsymbol{\Lambda}_{1}=\mathrm{Diag}\left(\mathbf{a}\right)$
being diagonal. The affine mapping from $\mathbf{x}$ to $\mathbf{y}$
is also analogous to Theorem \ref{thm:eqv recast}. Let $\mathbf{y}=\mathbf{T}^{T}\mathbf{x}+\boldsymbol{\Lambda}_{1}^{-1}\mathbf{T}^{-1}\mathbf{b}_{1}$,
$\mathbf{b}=\boldsymbol{\Lambda}_{1}^{-1}\mathbf{T}^{-1}\mathbf{b}_{1}-\mathbf{T}^{-1}\mathbf{b}_{0}$,
and $c=\left(\mathbf{T}^{-1}\mathbf{b}_{1}\right)^{T}\boldsymbol{\Lambda}_{1}^{-1}\left(\mathbf{T}^{-1}\mathbf{b}_{1}\right)-c_{1}$.
Then we show the mapping:
\begin{equation}
\begin{cases}
\mathbf{x}^{T}\mathbf{A}_{0}\mathbf{x}+2\mathbf{b}_{0}^{T}\mathbf{x}+c_{0}=\mathbf{y}^{T}\mathbf{y}-2\mathbf{b}^{T}\mathbf{y}+\mathrm{const}\\
\mathbf{x}^{T}\mathbf{A}_{1}\mathbf{x}+2\mathbf{b}_{1}^{T}\mathbf{x}+c_{1}=\mathbf{y}^{T}\boldsymbol{\Lambda}_{1}\mathbf{y}-c.
\end{cases}
\end{equation}
The original QCQP is decoupled as 
\begin{equation}
\begin{aligned} & \underset{\left\{ y_{i}\right\} _{i=1}^{n}}{\mathsf{minimize}} &  & \sum_{i=1}^{N}y_{i}^{2}-2b_{i}y_{i}\\
 & \mathsf{subject\thinspace to} &  & \sum_{i=1}^{N}a_{i}y_{i}^{2}-c=0.
\end{aligned}
\end{equation}

Now we are prepared for the bisection search stage. Note that $a_{i}$
could be positive or negative, and we can assert that $\max\left(\left\{ a_{i}\right\} \right)$
is the largest positive value, and that $\min\left(\left\{ a_{i}\right\} \right)$
is the smallest negative value. According to the first-order KKT conditions,
the optimal $y_{i}$'s are 
\begin{equation}
y_{i}=b_{i}/\left(1+\lambda a_{i}\right).
\end{equation}
The dual variable $\lambda$ satisfies $f\left(\lambda\right)=c$
where 
\begin{equation}
f\left(\lambda\right)=\sum_{i=1}^{N}\left(a_{i}b_{i}^{2}\right)/\left(1+\lambda a_{i}\right)^{2}.
\end{equation}
We further apply the second-order KKT condition to narrow the range
of $\lambda$: 
\begin{equation}
\begin{aligned} & 2+\lambda\cdot2a_{i}\geq0,\forall i\\
\Longrightarrow & -1/\max\left(\left\{ a_{i}\right\} \right)\leq\lambda\leq-1/\min\left(\left\{ a_{i}\right\} \right).
\end{aligned}
\end{equation}
Within $\left(-1/\max\left(\left\{ a_{i}\right\} \right),-1/\min\left(\left\{ a_{i}\right\} \right)\right)$,
$f\left(\lambda\right)$ is monotonically decreasing: 
\begin{equation}
f^{\prime}\left(\lambda\right)=-\sum_{i=1}^{N}2\left(a_{i}^{2}b_{i}^{2}\right)\left(1+\lambda a_{i}\right)^{-3}\leq0.
\end{equation}
Moreover, $f\left(\lambda\right)\in\mathbb{R}$ in the specified interval,
so the optimal $\lambda$ can still be solved from bisection search. 

To conclude, in the indefinite scenario, we solve QCQP \eqref{eq: QCQP eq real}
under the assumptions $\mathbf{A}_{0}\succ\mathbf{0}$, and $\mathbf{A}_{1}$
is symmetric, indefinite, and full-rank. The proposed solution scheme
is still applicable with minor modifications to the affine mapping
and the dual variable equation. 

\section{Extension III: Constraint Augmentation\label{sec:Constraint-Augmentation-Extension}}

In this section, we introduce additional linear equalities as an extension
to QCQP \eqref{eq: QCQP eq real}. The linear constraints are compactly
written as $\mathbf{A}_{2}\mathbf{x}=\mathbf{b}_{2}$ where $\mathbf{A}_{2}\in\mathbb{R}^{p\times N}$
is a rectangular matrix with $p\leq N$ and $\mathbf{b}_{2}\in\mathbb{R}^{p}$.
In engineering practice, linear equalities are independent and thus
$\mathbf{A}_{2}$ has full row rank. We inherit the assumptions in
Sec. \ref{sec:Globally-Optimal-Solution} on $\mathbf{A}_{0}$ and
$\mathbf{A}_{1}$: $\mathbf{A}_{0}$ is symmetric and $\mathbf{A}_{1}\succ\mathbf{0}$.
After constraint augmentation, the QCQP problem becomes 
\begin{equation}
\begin{aligned} & \underset{\mathbf{x}\in\mathbb{R}^{N}}{\mathsf{minimize}} &  & \mathbf{x}^{T}\mathbf{A}_{0}\mathbf{x}+2\mathbf{b}_{0}^{T}\mathbf{x}+c_{0}\\
 & \mathsf{subject\thinspace to} &  & \mathbf{x}^{T}\mathbf{A}_{1}\mathbf{x}+2\mathbf{b}_{1}^{T}\mathbf{x}+c_{1}=0\\
 &  &  & \mathbf{A}_{2}\mathbf{x}=\mathbf{b}_{2}.
\end{aligned}
\label{eq:QCQP eq real ext1}
\end{equation}
QCQP \eqref{eq:QCQP eq real ext1} has a simple prototype in \citep[eq. (21)]{zhao2018mean}. 

It is generally acknowledged that the solution to any feasible linear
system can be expressed as the sum of a particular solution and a
homogeneous one. Since $\mathbf{A}_{2}$ has full row rank, we can
easily find a particular solution $\mathbf{x}_{p}$:
\begin{equation}
\mathbf{x}_{p}=\mathbf{A}_{2}^{T}\left(\mathbf{A}_{2}\mathbf{A}_{2}^{T}\right)^{-1}\mathbf{b}_{2}.
\end{equation}
The homogeneous solution $\mathbf{x}_{h}$ can be constructed from
the null space of $\mathbf{A}_{2}$, represented by $\mathbf{A}_{2\perp}$
($\mathbf{A}_{2}\mathbf{A}_{2\perp}=\mathbf{0}$ and $\mathbf{A}_{2\perp}$
has full column rank): 
\begin{equation}
\mathbf{x}_{h}=\mathbf{A}_{2\perp}\tilde{\mathbf{x}},\tilde{\mathbf{x}}\in\mathbb{R}^{N-p}.
\end{equation}
Therefore, we can eliminate the affine constraints $\mathbf{A}_{2}\mathbf{x}=\mathbf{b}_{2}$
using $\mathbf{x}=\mathbf{x}_{p}+\mathbf{A}_{2\perp}\tilde{\mathbf{x}}$
and transform \eqref{eq:QCQP eq real ext1} back to a single-equality-constraint
QCQP with variable $\tilde{\mathbf{x}}$: 
\begin{equation}
\begin{aligned} & \underset{\tilde{\mathbf{x}}\in\mathbb{R}^{N-p}}{\mathsf{minimize}} &  & \tilde{\mathbf{x}}^{T}\tilde{\mathbf{A}}_{0}\tilde{\mathbf{x}}+2\tilde{\mathbf{b}}_{0}^{T}\tilde{\mathbf{x}}+\tilde{c}_{0}\\
 & \mathsf{subject\thinspace to} &  & \tilde{\mathbf{x}}^{T}\tilde{\mathbf{A}}_{1}\tilde{\mathbf{x}}+2\tilde{\mathbf{b}}_{1}^{T}\tilde{\mathbf{x}}+\tilde{c}_{1}=0,
\end{aligned}
\label{eq:QCQP ext1 reduce}
\end{equation}
where $\tilde{\mathbf{A}}_{0}=\mathbf{A}_{2\perp}^{T}\mathbf{A}_{0}\mathbf{A}_{2\perp}$,
$\tilde{\mathbf{b}}_{0}=\mathbf{A}_{2\perp}^{T}\left(\mathbf{A}_{0}\mathbf{x}_{p}+\mathbf{b}_{0}\right)$,
$\tilde{c}_{0}=\mathbf{x}_{p}^{T}\mathbf{A}_{0}\mathbf{x}_{p}+2\mathbf{b}_{0}^{T}\mathbf{x}_{p}+c_{0}$,
$\tilde{\mathbf{A}}_{1}=\mathbf{A}_{2\perp}^{T}\mathbf{A}_{1}\mathbf{A}_{2\perp}$,
$\tilde{\mathbf{b}}_{1}=\mathbf{A}_{2\perp}^{T}\left(\mathbf{A}_{1}\mathbf{x}_{p}+\mathbf{b}_{1}\right)$,
and $\tilde{c}_{1}=\mathbf{x}_{p}^{T}\mathbf{A}_{1}\mathbf{x}_{p}+2\mathbf{b}_{1}^{T}\mathbf{x}_{p}+c_{1}$.
It is straightforward to see that a symmetric $\mathbf{A}_{0}$ yields
a symmetric $\tilde{\mathbf{A}}_{0}$, and that $\mathbf{A}_{1}\succ\mathbf{0}$
and $\mathbf{A}_{2\perp}$ being of full column rank imply that $\tilde{\mathbf{A}}_{1}\succ\mathbf{0}$.
If $-\tilde{\mathbf{b}}_{1}^{T}\tilde{\mathbf{A}}_{1}^{-1}\tilde{\mathbf{b}}_{1}+\tilde{c}_{1}\leq0$,
the constraint set is nonempty and QCQP \eqref{eq:QCQP ext1 reduce}
fits in the ``SD-BS'' solution scheme.

\section{Extension IV: Matrix-Format Variable\label{sec:Matrix-Format-Extension}}

In wireless communications, QCQP \eqref{eq: QCQP eq real} may be
expressed in a complex form, with the optimization variable in matrix
format: 

\begin{equation}
\begin{aligned} & \underset{\mathbf{X}\in\mathbb{C}^{N_{1}\times N_{2}}}{\mathsf{minimize}} &  & \mathrm{Tr}\left(\mathbf{X}^{H}\mathbf{A}_{0}\mathbf{X}\right)+2\mathrm{Re}\left[\mathrm{Tr}\left(\mathbf{B}_{0}^{H}\mathbf{X}\right)\right]+c_{0}\\
 & \mathsf{subject\thinspace to} &  & \mathrm{Tr}\left(\mathbf{X}^{H}\mathbf{A}_{1}\mathbf{X}\right)+2\mathrm{Re}\left[\mathrm{Tr}\left(\mathbf{B}_{1}^{H}\mathbf{X}\right)\right]+c_{1}=0,
\end{aligned}
\label{eq:complex and matrix format}
\end{equation}
with $\mathbf{A}_{0}=\mathbf{A}_{0}^{H}$ and $\mathbf{A}_{1}\succ\mathbf{0}$.
A conventional way to solve \eqref{eq:complex and matrix format},
as has been revealed in \citep[eq. (35)]{yu2016alternating}, is to
apply the vectorization technique on the variable $\mathbf{X}$. The
overall complexity is $\mathcal{O}\left(\left(N_{1}N_{2}\right)^{6.5}\right)$
due to SDR, which is extremely time-consuming for solving a subproblem.
But if we choose the ``SD-BS'' scheme instead, the order of complexity
can be reduced to $\mathcal{O}\left(N_{1}^{3}\right)$. We will explain
how to achieve this improvement.

The first stage is very similar to the standard case. We perform simultaneous
diagonalization in the same way as Lemma \ref{lem:SD-full}: $\mathbf{A}_{0}=\mathbf{T}\boldsymbol{\Lambda}_{0}\mathbf{T}^{H}$
and $\mathbf{A}_{1}=\mathbf{T}\mathbf{T}^{H}$. Then we derive an
affine mapping $\mathbf{Y}=\mathbf{T}^{H}\mathbf{X}+\mathbf{T}^{-1}\mathbf{B}_{1}$,
analogous to the one in Theorem \ref{thm:eqv recast}, and transform
the constraint and objective as 
\begin{equation}
\mathrm{Tr}\left(\mathbf{X}^{H}\mathbf{A}_{1}\mathbf{X}\right)+2\mathrm{Re}\left[\mathrm{Tr}\left(\mathbf{B}_{1}^{H}\mathbf{X}\right)\right]+c_{1}=\mathrm{Tr}\left(\mathbf{Y}^{H}\mathbf{Y}\right)-c
\end{equation}
and 
\begin{equation}
\begin{aligned} & \mathrm{Tr}\left(\mathbf{X}^{H}\mathbf{A}_{0}\mathbf{X}\right)+2\mathrm{Re}\left[\mathrm{Tr}\left(\mathbf{B}_{0}^{H}\mathbf{X}\right)\right]+c_{0}\\
= & \mathrm{Tr}\left(\mathbf{Y}^{H}\boldsymbol{\Lambda}_{0}\mathbf{Y}\right)-2\mathrm{Re}\left[\mathrm{Tr}\left(\mathbf{B}^{H}\mathbf{Y}\right)\right]+\mathrm{const}
\end{aligned}
\end{equation}
with $c=\mathrm{Tr}\left(\left(\mathbf{T}^{-1}\mathbf{B}_{1}\right)^{H}\left(\mathbf{T}^{-1}\mathbf{B}_{1}\right)\right)-c_{1}$,
$\mathbf{B}=\boldsymbol{\Lambda}_{0}\mathbf{T}^{-1}\mathbf{B}_{1}-\mathbf{T}^{-1}\mathbf{B}_{0}$,
and $\mathrm{const}$ is an optimization irrelevant constant. Hence,
the matrix-format QCQP is recast as
\begin{equation}
\begin{aligned} & \underset{\mathbf{Y}}{\mathsf{minimize}} &  & \mathrm{Tr}\left(\mathbf{Y}^{H}\boldsymbol{\Lambda}_{0}\mathbf{Y}\right)-2\mathrm{Re}\left[\mathrm{Tr}\left(\mathbf{B}^{H}\mathbf{Y}\right)\right]\\
 & \mathsf{subject\thinspace to} &  & \mathrm{Tr}\left(\mathbf{Y}^{H}\mathbf{Y}\right)-c=0,
\end{aligned}
\label{eq:decoupled QCQP}
\end{equation}
or the equivalent decoupled formulation 
\begin{equation}
\begin{aligned} & \underset{\left\{ \mathbf{y}_{i}\right\} _{i=1}^{N_{2}}}{\mathsf{minimize}} &  & \sum_{i=1}^{N_{2}}\mathbf{a}^{T}\left|\mathbf{y}_{i}\right|^{2}-2\mathrm{Re}\left[\mathbf{b}_{i}^{H}\mathbf{y}_{i}\right]\\
 & \mathsf{subject\thinspace to} &  & \sum_{i=1}^{N_{2}}\left\Vert \mathbf{y}_{i}\right\Vert _{2}^{2}-c=0,
\end{aligned}
\label{eq:decoupled QCQP summation}
\end{equation}
where we undo $\boldsymbol{\Lambda}_{0}$ as $\mathrm{Diag}\left(\mathbf{a}\right)$
($\mathbf{a}\in\mathbb{R}^{N_{1}}$) and split $\mathbf{Y}$ and $\mathbf{B}$
columnwisely as $\mathbf{Y}=\left[\mathbf{y}_{1},\ldots,\mathbf{y}_{N_{2}}\right]$
and $\mathbf{B}=\left[\mathbf{b}_{1},\ldots,\mathbf{b}_{N_{2}}\right]$. 

The second stage begins with an examination of the KKT conditions.
The first-order KKT solution to the decoupled QCQP \eqref{eq:decoupled QCQP summation}
is presented in the theorem below. 
\begin{thm}
\label{thm:KKT mat}The KKT solution to QCQP \eqref{eq:decoupled QCQP summation},
denoted by $\tilde{\mathbf{Y}}$, is expressed columnwisely as: $\forall i$,
\begin{equation}
\arg\left(\tilde{\mathbf{y}}_{i}\right)=\arg\left(\mathbf{b}_{i}\right)
\end{equation}
and 
\begin{equation}
\left|\tilde{\mathbf{y}}_{i}\right|=\left|\mathbf{b}_{i}\right|/\left(\mathbf{a}+\lambda\mathbf{1}\right)
\end{equation}
where $\lambda$ satisfies $f\left(\lambda\right)=\sum_{j=1}^{N_{1}}b_{j}^{2}/\left(a_{j}+\lambda\right)^{2}$
and $b_{j}=\sqrt{\sum_{i=1}^{N_{2}}\left|b_{i,j}\right|^{2}}$. 
\end{thm}
\begin{IEEEproof}
See Appendix \ref{sec:Proof-of-Theorem KKT mat} for the detailed
proof. 
\end{IEEEproof}
We find that $f\left(\lambda\right)$ in Theorem \ref{thm:KKT mat}
has exactly the same form as \eqref{eq:f lambda}. As a result, the
optimal $\lambda$ can be solved in the same way as described in Section
\ref{subsec:Bisection-Search}, using the second-order KKT condition. 

Specially, \citet{yu2016alternating} presented an even simpler example
of single-equality-constraint matrix-format QCQP: 
\begin{equation}
\begin{aligned} & \underset{\mathbf{X}\in\mathbb{C}^{N_{1}\times N_{2}}}{\mathsf{minimize}} &  & \left\Vert \mathbf{F}-\mathbf{G}\mathbf{X}\right\Vert _{F}^{2}\\
 & \mathsf{subject\thinspace to} &  & \left\Vert \mathbf{X}\right\Vert _{F}^{2}=c.
\end{aligned}
\label{eq:QCQP example}
\end{equation}
In this example, $\mathbf{A}_{0}=\mathbf{G}^{H}\mathbf{G}\succ\mathbf{0}$,
$\mathbf{B}_{0}=\mathbf{G}^{H}\mathbf{F}$, $\mathbf{A}_{1}=\mathbf{I}$,
and \textbf{$\mathbf{B}_{1}=\mathbf{0}$}. Simultaneous diagonalization
can be simplified with an SVD operation on $\mathbf{G}$. With $\mathbf{G}=\mathbf{U}\boldsymbol{\Sigma}\mathbf{V}^{H}$,
we obtain ($\boldsymbol{\Lambda}_{0}=\boldsymbol{\Sigma}^{H}\boldsymbol{\Sigma}$,
still diagonal)
\begin{equation}
\begin{cases}
\mathbf{A}_{0}=\mathbf{V}\boldsymbol{\Lambda}_{0}\mathbf{V}^{H}\\
\mathbf{A}_{1}=\mathbf{V}\mathbf{V}^{H},
\end{cases}
\end{equation}
so $\mathbf{T}=\mathbf{V}$, which is unitary. Given all these parameters,
we further have 
\begin{equation}
\mathbf{Y}=\mathbf{T}^{H}\mathbf{X}+\mathbf{T}^{-1}\mathbf{B}_{1}=\mathbf{V}^{H}\mathbf{X}
\end{equation}
and 
\begin{equation}
\mathbf{B}=\boldsymbol{\Lambda}_{0}\mathbf{T}^{-1}\mathbf{B}_{1}-\mathbf{T}^{-1}\mathbf{B}_{0}=-\mathbf{V}^{H}\mathbf{G}^{H}\mathbf{F}.
\end{equation}
After reformulation, QCQP \eqref{eq:QCQP example} becomes exactly
\eqref{eq:decoupled QCQP} with $\boldsymbol{\Lambda}_{0}\succeq\mathbf{0}$
(i.e., $\mathbf{a}$ is elementwisely nonnegative). The rest of the
solution has already been elaborated on above. 

\section{Numerical Simulations \label{sec:Numerical-Simulations}}

In this section, we present numerical results on solving several classes
of single-equality-constraint QCQP problems. All simulations are
performed on a PC with a 2.90 GHz i7-10700 CPU and 16.0 GB RAM. We
will perform a detailed comparison between the proposed approach and
various existing methods. A prevalent benchmark method is SDR because
of tight relaxation in the case of single-equality-constraint QCQPs
and SDR is implemented with SDP solvers whose algorithm core is the
interior point method. In our simulation, we use the solver ``MOSEK''
for SDR realization. Apart from the SDR, we also include a built-in
solver of Matlab, ``fmincon'', as a general benchmark for evaluation. 

The compared algorithms will be measured in terms of constraint satisfaction,
optimality gap, and computational time. Constraint satisfaction is
defined as 
\begin{equation}
\textrm{Constraint Satisfaction}=\left|\textrm{Constraint}\left(\mathbf{x}_{\textrm{out}}\right)\right|.
\end{equation}
Optimality gap is calculated as 
\begin{equation}
\begin{aligned} & \textrm{Optimality Gap}=\\
 & \left|\textrm{Objective}\left(\mathbf{x}_{\textrm{out}}\right)-\textrm{duality objective bound}\right|.
\end{aligned}
\end{equation}
The duality objective bound is the optimal value of the dual problem
and can be retrieved from the outputs of an off-the-shelf solver if
the problem is properly solved with sufficient accuracy. The reported
performance records are all averaged from $100$ randomized instances.
If a data point is absent, the particular algorithm fails to produce
valid results for the given setting. 
\begin{figure*}[tbh]
\begin{centering}
\subfloat[\label{fig:Constraint-satisfaction-1}Constraint satisfaction.]{\begin{centering}
\includegraphics[width=0.205\paperwidth]{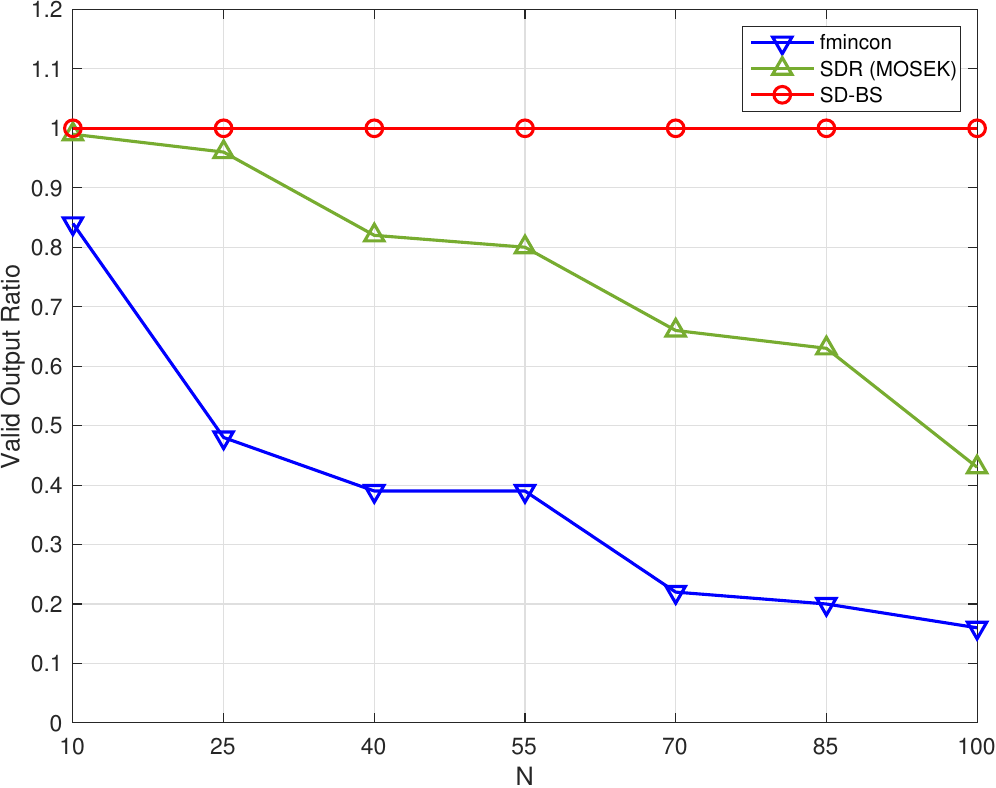}
\par\end{centering}
}\hspace{0.05\paperwidth}\subfloat[\label{fig:Optimality-gap-1}Optimality gap.]{\begin{centering}
\includegraphics[width=0.205\paperwidth]{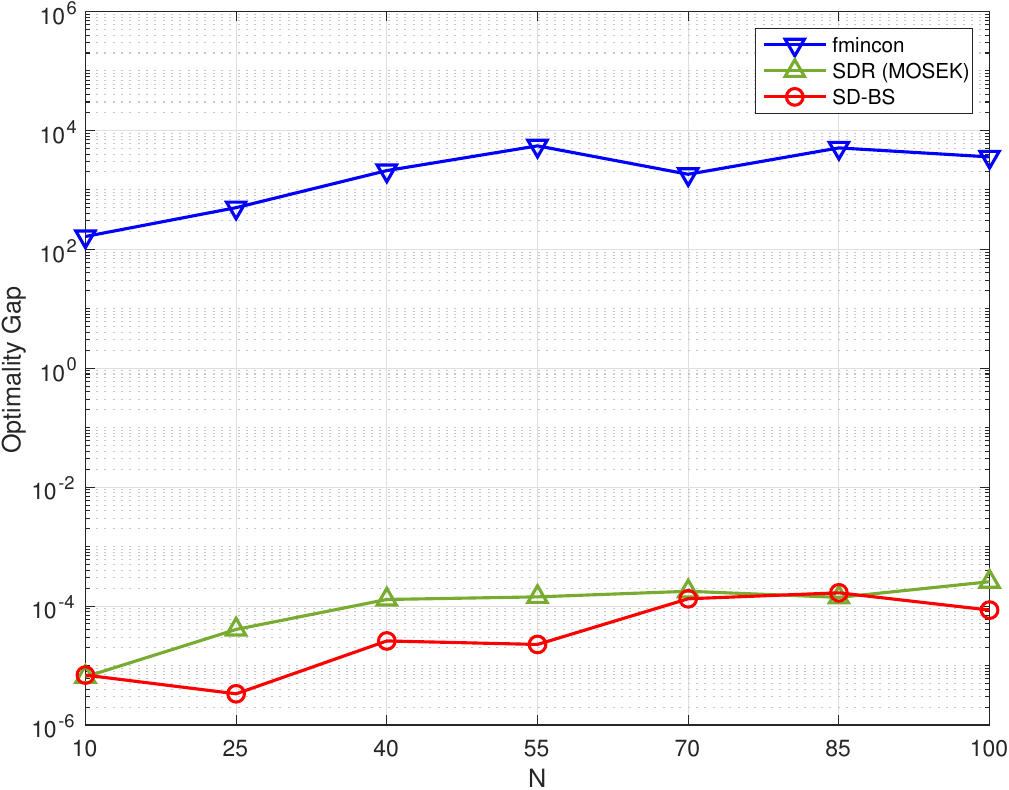}
\par\end{centering}
}\hspace{0.05\paperwidth}\subfloat[\label{fig:Computational-time-1}Computational time.]{\begin{centering}
\includegraphics[width=0.205\paperwidth]{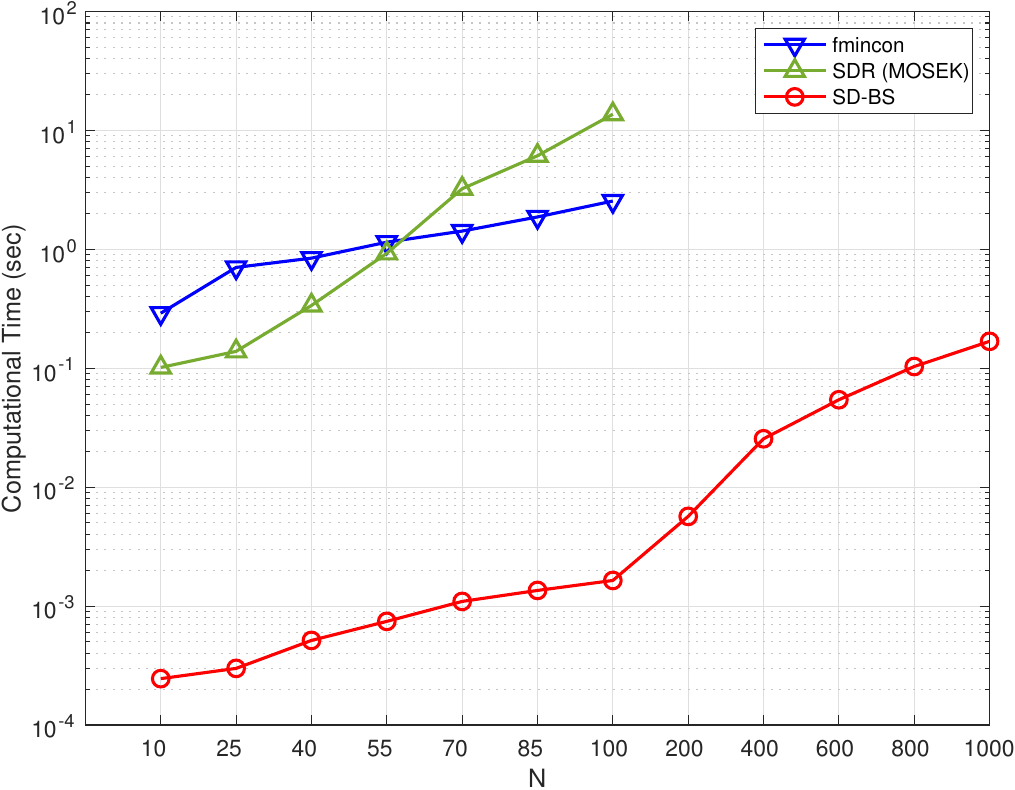}
\par\end{centering}
}
\par\end{centering}
\caption{\label{fig:main_comp_1_opt_time}Performance metrics versus problem
dimension $N$. $\mathbf{A}_{0}$ is real symmetric and $\mathbf{A}_{1}\succ\mathbf{0}$.}
\end{figure*}
\begin{figure*}[t]
\begin{centering}
\subfloat[\label{fig:Constraint-satisfaction-2}Constraint satisfaction.]{\begin{centering}
\includegraphics[width=0.205\paperwidth]{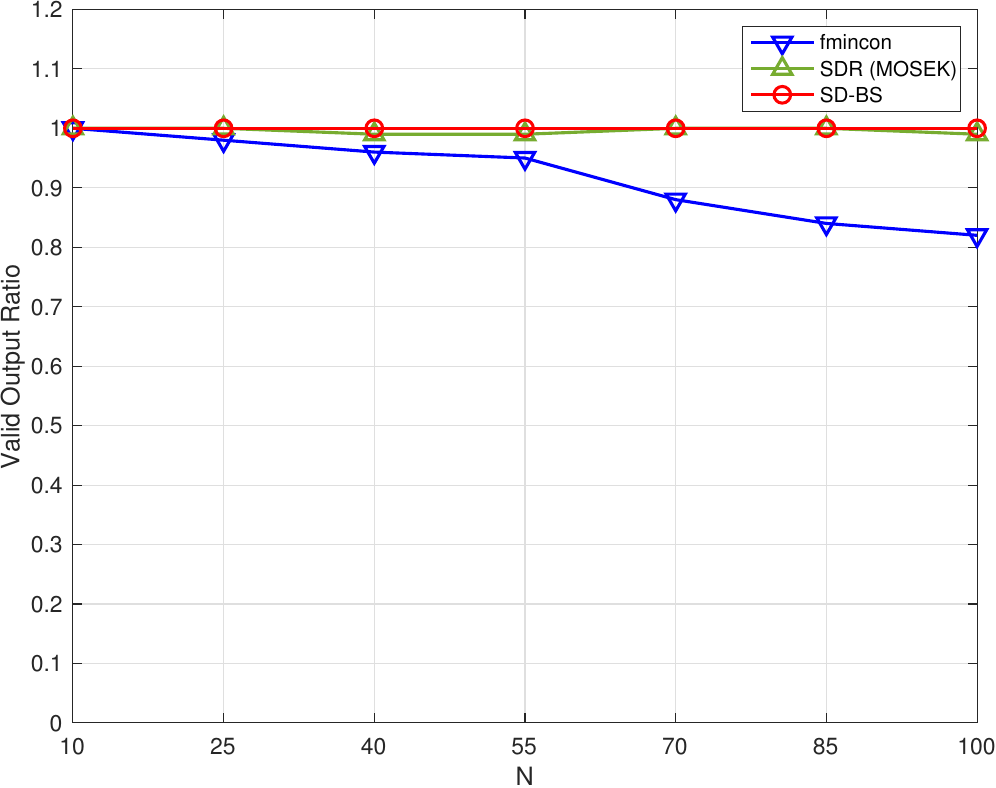}
\par\end{centering}
}\hspace{0.05\paperwidth}\subfloat[\label{fig:Optimality-gap-2}Optimality gap.]{\begin{centering}
\includegraphics[width=0.205\paperwidth]{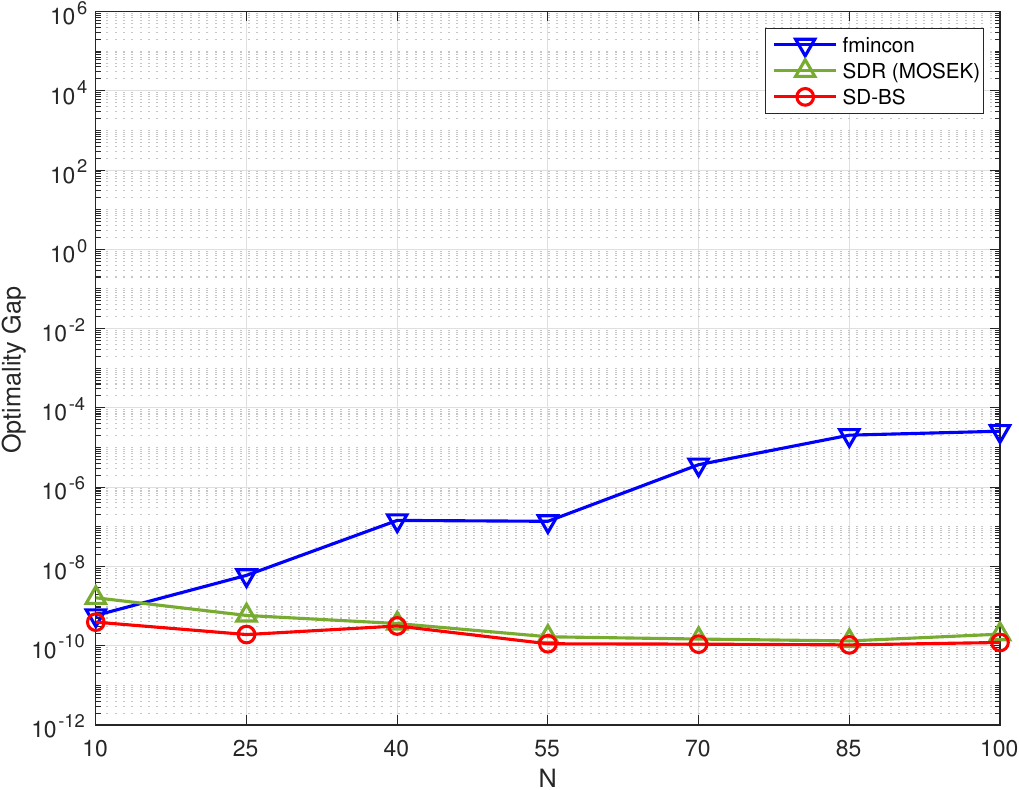}
\par\end{centering}
}\hspace{0.05\paperwidth}\subfloat[\label{fig:Computational-time-2}Computational time.]{\begin{centering}
\includegraphics[width=0.205\paperwidth]{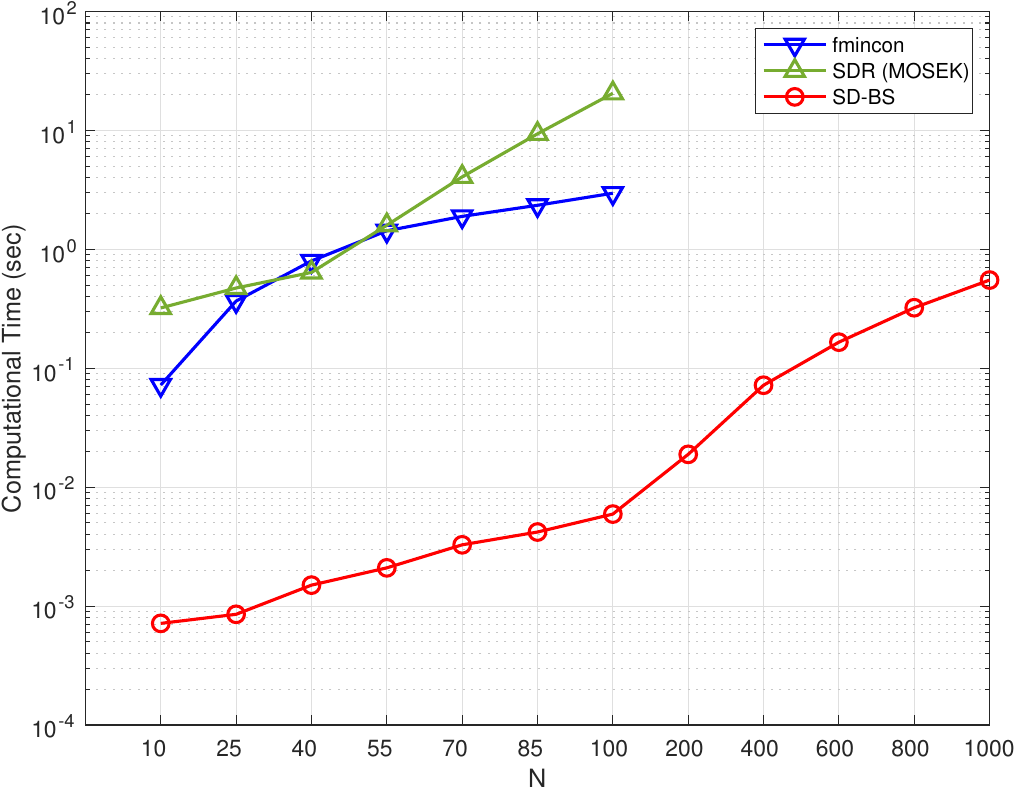}
\par\end{centering}
}
\par\end{centering}
\caption{\label{fig:main_comp_2_opt_time}Performance metrics versus problem
dimension $N$. $\mathbf{A}_{0}\succeq\mathbf{0}$, $\mathbf{A}_{1}\succeq\mathbf{0}$
and $\mathbf{A}_{1}\nsucc\mathbf{0}$.}
\end{figure*}

\subsection{Standard Single-Equality-Constraint QCQP \label{subsec:Standard-Single-Equality-QCQP}}

We generate standard QCQP instances in the following way. The positive
definite parameter $\mathbf{A}_{1}\in\mathbb{R}^{N\times N}$ is generated
as $\tau_{1}\mathbf{F}_{1}\mathbf{F}_{1}^{T}+\varepsilon\mathbf{I}$
where $\tau_{1}\sim\chi^{2}$, $\mathbf{F}_{1}\in\mathbb{R}^{N\times N}$
is a multiplicative factor with independent and identically distributed
(i.i.d.) Gaussian entries drawn from $\mathcal{N}\left(0,1\right)$,
and $\varepsilon$ is a positive constant with a default choice of
$0.1$. $\mathbf{A}_{0}\in\mathbb{R}^{N\times N}$ is real symmetric:
$\mathbf{A}_{0}=\tau_{0}\mathbf{F}_{0}\textrm{Diag}\left(\mathbf{f}_{0}\right)\mathbf{F}_{0}^{T}$
where $\tau_{0}\sim\chi^{2}$, the elements of $\mathbf{F}_{0}\in\mathbb{R}^{N\times N}$
are also i.i.d. Gaussian distributed, each drawn from $\mathcal{N}\left(0,1\right)$,
and vector entries in$\mathbf{f}_{0}\in\mathbb{R}^{N\times1}$ are
uniformly distributed between $-1$ and $1$. The entries of $\mathbf{b}_{0}\in\mathbb{R}^{N\times1}$
and $\mathbf{b}_{1}\in\mathbb{R}^{N\times1}$ are i.i.d. random variables,
each following $\mathcal{N}\left(0,1\right)$. $c_{0}$ is optimization
invariant and set to zero without loss of generality. $c_{1}$ is
chosen such that $c_{1}\leq\mathbf{b}_{1}^{T}\mathbf{A}_{1}^{-1}\mathbf{b}_{1}$
for the sake of solution existence. The compared algorithms include:
SDR (MOSEK), fmincon, SD-BS (proposed). 

When $\mathbf{A}_{0}$ is real symmetric and $\mathbf{A}_{1}\succ\mathbf{0}$,
the comparison results are presented in Figure \ref{fig:main_comp_1_opt_time}.
First, we assess the performance on constraint satisfaction. Valid
outputs are defined as those whose constraint satisfaction level is
no larger than $10^{-5}$. Figure \ref{fig:Constraint-satisfaction-1}
displays the ratio of valid outputs for the compared algorithms. Of
the two benchmarks (SDR and fmincon), the validity ratio shows a descending
trend as the problem size $N$ increases. When $N$ reaches 100, only
$40\%$ and $20\%$ of the QCQP solutions strictly meet the equality
constraint for SDR and fmincon. However, all the solutions produced
by SD-BS are valid irrespective of the choice of $N$. 

Among the valid solutions, we further examine the metrics of optimality
gap and computation time. In Figure \ref{fig:Optimality-gap-1}, we
can see that the general Matlab solver fmincon suffers a huge optimality
gap, about six orders of magnitude larger than the rest. In contrast,
the optimality gaps of SDR and SD-BS are well controlled, remaining
mostly within $10^{-4}$. The performance of computational efficiency
is demonstrated in Figure \ref{fig:Computational-time-1}. Due to
the limited capability of solving large-scale problems, the problem
sizes for fmincon and SDR are restricted to $N\leq100$. The proposed
algorithm SD-BS is at least two orders of magnitude faster than the
other benchmarks when $N\leq100$. Moreover, SD-BS can support problem
sizes up to $10^{3}$, with an average running time of less than $0.2$
seconds. These two observations indicate better computational efficiency
and greater scalability (at least ten times) compared to the existing
benchmarks. 

\subsection{Rank Deficiency Extension}

In the rank deficient scenario, we generate $\mathbf{A}_{1}$ as $\tau_{1}\mathbf{F}_{1}\mathbf{F}_{1}^{T}$
and $\mathbf{A}_{0}$ as $\tau_{0}\mathbf{F}_{0}\textrm{Diag}\left(\mathbf{f}_{0}\right)\mathbf{F}_{0}^{T}$
where the generation of $\tau_{1}$, $\mathbf{F}_{1}\in\mathbb{R}^{N\times r}$,
$\tau_{0}$ and $\mathbf{F}_{0}\in\mathbb{R}^{N\times N}$ follows
Section \ref{subsec:Standard-Single-Equality-QCQP} and vector entries
in $\mathbf{f}_{0}$ are uniformly distributed between $0$ and $1$.
We choose the rank of $\mathbf{A}_{1}$ as $r=\mathrm{floor}\left(N/2\right)$. 

In terms of constraint satisfaction, as shown in Figure \ref{fig:Constraint-satisfaction-2},
the outputs of SDR and SD-BS are almost $100\%$ valid across the
entire range of $N$ while the validity ratio of fmincon still exhibits
a downward trend. When $N=100$, around $80\%$ of the solutions strictly
meet the feasibility constraint. The metrics of optimality gap and
computation time are measured across the valid solutions. Figure \ref{fig:Optimality-gap-2}
implies that of the compared algorithms, all the valid solutions are
approximately optimal in the numerical sense. The optimality gaps
of SDR and SD-BS are extremely close to zero (lower than $10^{-8}$)
but fmincon performs slightly inferior, reaching a gap level of $10^{-4}$
at $N=100$. Since there is no big difference in the optimality metric,
computational efficiency becomes the top concern. As can be observed
in Figure \ref{fig:Computational-time-2}, the proposed algorithm
SD-BS is about two orders of magnitude faster than the benchmarks
and offers a tenfold wider range of scalability. The superiority of
SD-BS is thus justified. 

\subsection{Indefiniteness Extension}

\begin{figure*}[tbh]
\begin{centering}
\subfloat[\label{fig:Constraint-satisfaction-3}Constraint satisfaction.]{\begin{centering}
\includegraphics[width=0.205\paperwidth]{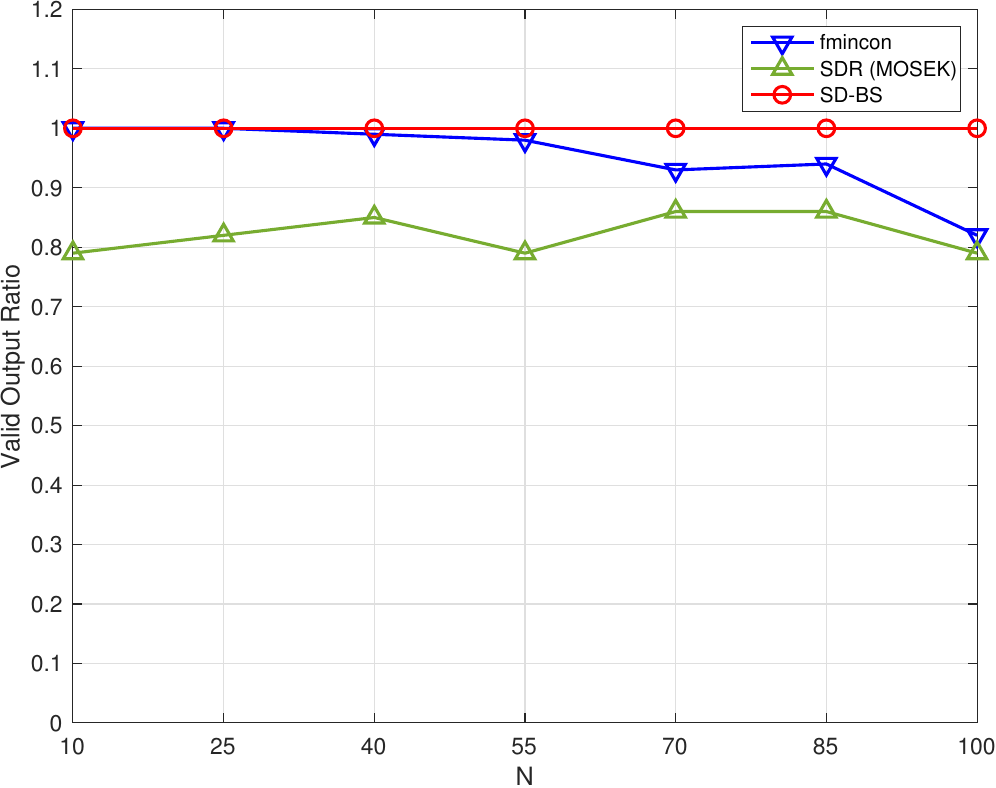}
\par\end{centering}
}\hspace{0.05\paperwidth}\subfloat[\label{fig:Optimality-gap-3}Optimality gap.]{\begin{centering}
\includegraphics[width=0.205\paperwidth]{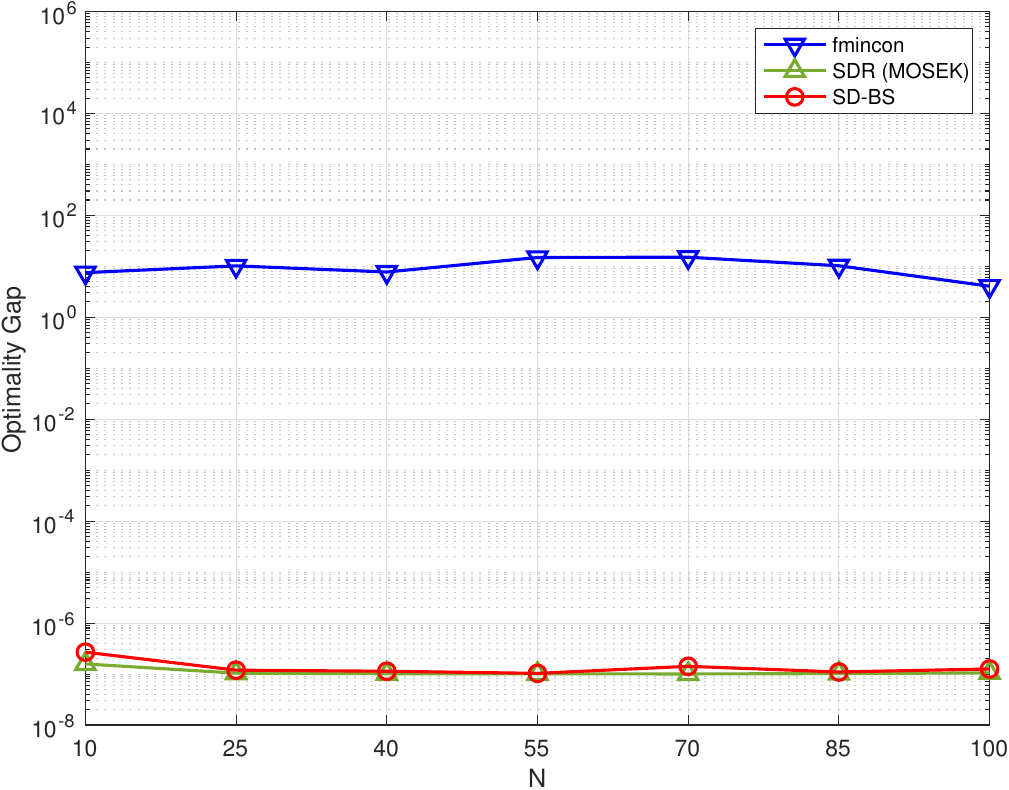}
\par\end{centering}
}\hspace{0.05\paperwidth}\subfloat[\label{fig:Computational-time-3}Computational time.]{\begin{centering}
\includegraphics[width=0.205\paperwidth]{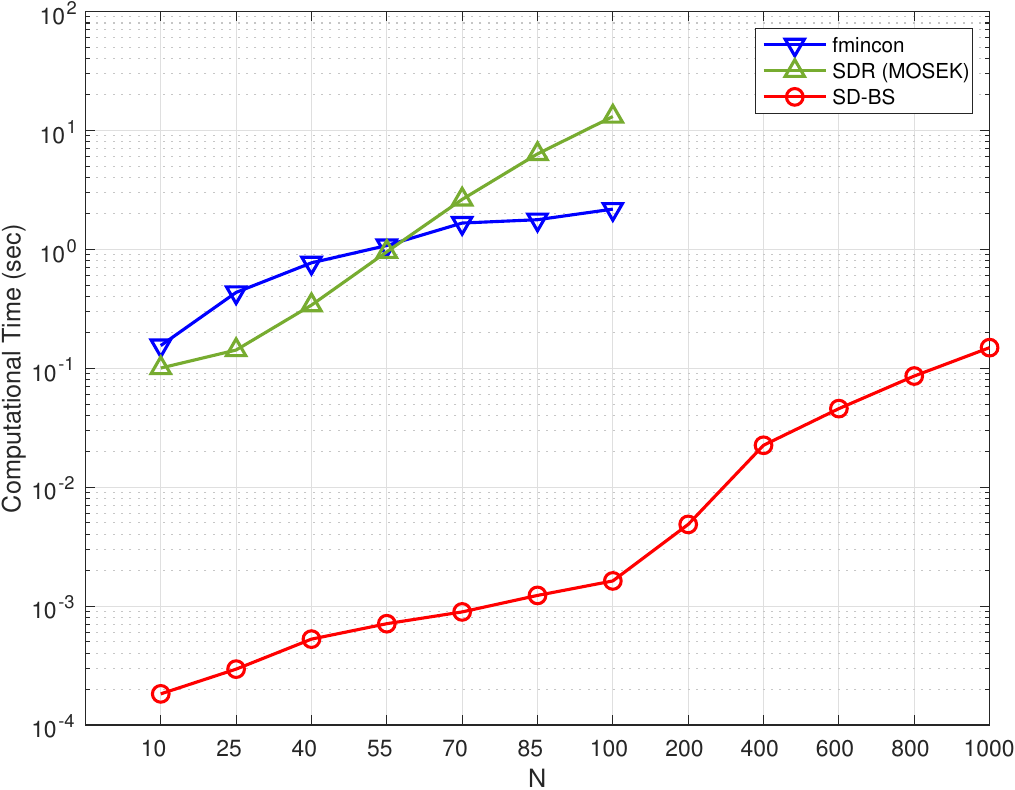}
\par\end{centering}
}
\par\end{centering}
\caption{\label{fig:main_comp_3_opt_time}Performance metrics versus problem
dimension $N$. $\mathbf{A}_{0}\succ\mathbf{0}$ and $\mathbf{A}_{1}$
is symmetric, indefinite, and full-rank.}
\end{figure*}
\begin{figure*}[tbh]
\begin{centering}
\subfloat[\label{fig:Constraint-satisfaction-4}Constraint satisfaction.]{\begin{centering}
\includegraphics[width=0.205\paperwidth]{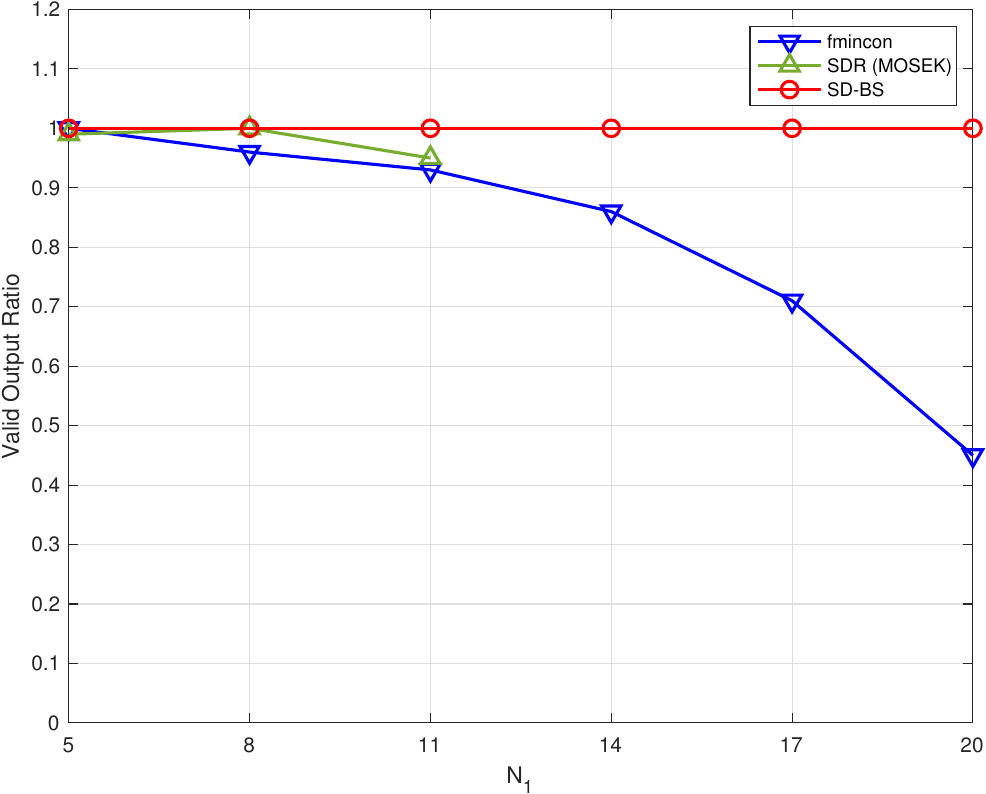}
\par\end{centering}
}\hspace{0.05\paperwidth}\subfloat[\label{fig:Optimality-gap-4}Optimality gap.]{\begin{centering}
\includegraphics[width=0.205\paperwidth]{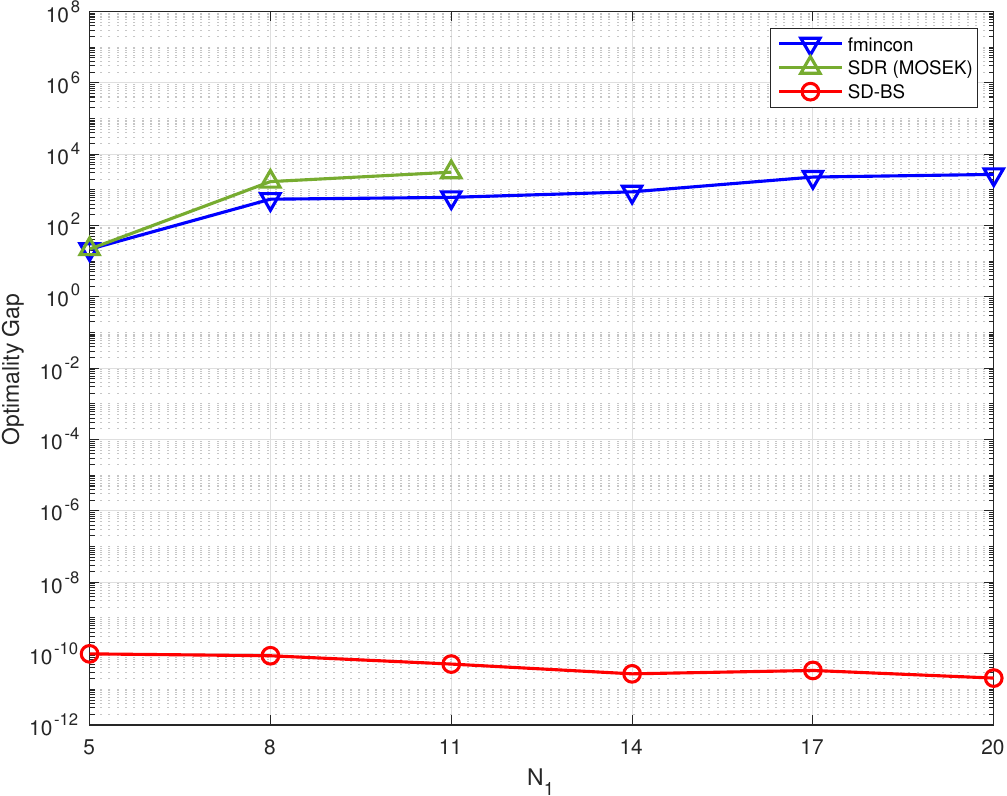}
\par\end{centering}
}\hspace{0.05\paperwidth}\subfloat[\label{fig:Computational-time-4}Computational time.]{\begin{centering}
\includegraphics[width=0.205\paperwidth]{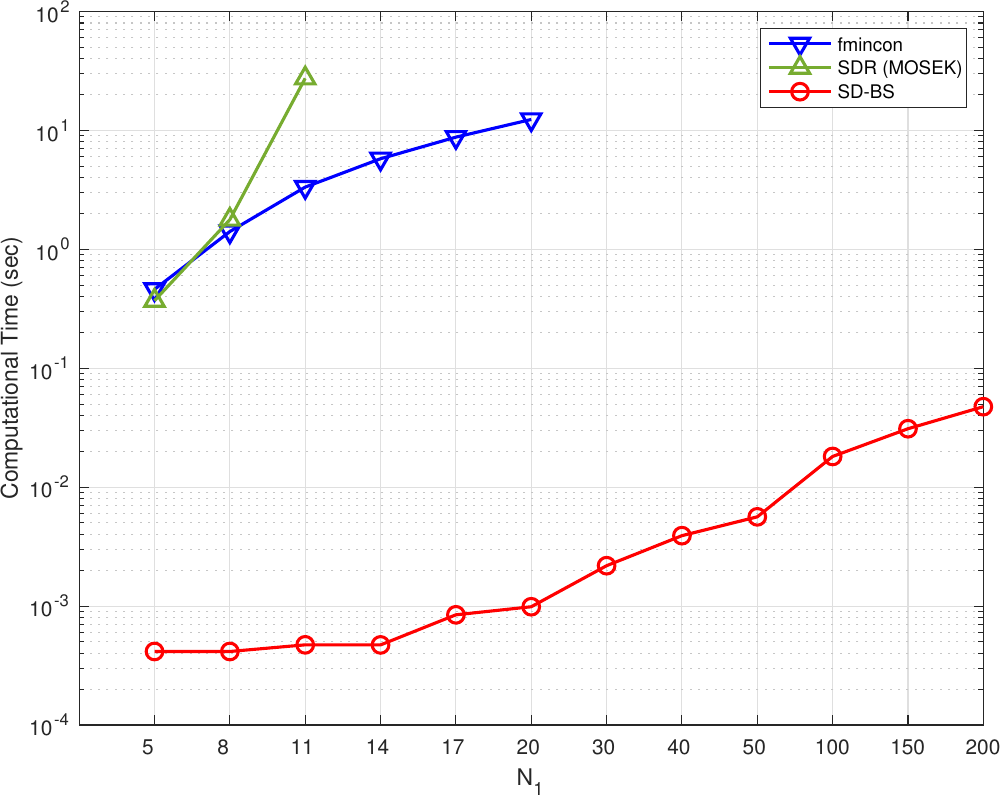}
\par\end{centering}
}
\par\end{centering}
\caption{\label{fig:main_comp_4_opt_time}Performance metrics versus problem
dimension $N_{1}$ for complex-valued matrix-format QCQP. $\mathbf{A}_{0}$
is Hermitian and $\mathbf{A}_{1}\succ\mathbf{0}$.}
\end{figure*}
In this scenario, $\mathbf{A}_{1}$ is symmetric, indefinite, and
full-rank, and $\mathbf{A}_{0}$ is positive definite. Thus, we generate
$\mathbf{A}_{0}$ as $\tau_{0}\mathbf{F}_{0}\mathbf{F}_{0}^{T}$ and
$\mathbf{A}_{1}$ as $\tau_{1}\mathbf{F}_{1}\textrm{Diag}\left(\mathbf{f}_{1}\right)\mathbf{F}_{1}^{T}$
where the generation of $\tau_{0}$, $\mathbf{F}_{0}\in\mathbb{R}^{N\times r}$,
$\tau_{1}$ and $\mathbf{F}_{1}\in\mathbb{R}^{N\times N}$ follows
Section \ref{subsec:Standard-Single-Equality-QCQP} and the elements
in $\mathbf{f}_{1}$ are uniformly distributed between $-1$ and $1$.
To ensure the full-rank property of $\mathbf{A}_{1}$, any $\mathbf{f}_{1}$
component, whose absolute value is smaller than $10^{-4}$, is set
to $10^{-4}$ while the original sign is kept. 

Figure \ref{fig:Constraint-satisfaction-3} describes the ratio of
valid outputs for different compared algorithms. The validity ratio
of SDR remains steady with the change in problem size and maintains
an $80\%$ validity rate. The valid solutions of fmincon slightly
decreases with $N$ and the validity rate is also over $80\%$. Every
solution of SD-BS satisfies the feasibility restriction under all
choices of $N$. In terms of optimality gap, the experimental phenomenon
in Figure \ref{fig:Optimality-gap-3} highly resembles Figure \ref{fig:Optimality-gap-1}:
the gap levels of SDR and SD-BS are as low as $10^{-7}$, whereas
the solutions provided by fmincon are clearly suboptimal. The numerical
results on computation time are displayed in Figure \ref{fig:Computational-time-3}.
The proposed SD-BS is at least two orders of magnitude more efficient
than the compared baselines. Moreover, SD-BS can solve large-scale
QCQPs with a size of $1000$, whereas the benchmarks can only handle
problems with $N\leq100$. 

\subsection{Matrix-Format Variable Extension}

Finally, we investigate into the matrix-format variable extension.
Recall that $\mathbf{X}\in\mathbb{C}^{N_{1}\times N_{2}}$ and we
set $N_{2}=N_{1}-1$. The parameters $\mathbf{A}_{0}\in\mathbb{C}^{N_{1}\times N_{1}}$,
$\mathbf{A}_{1}\in\mathbb{C}^{N_{1}\times N_{1}}$, $\mathbf{B}_{0}\in\mathbb{C}^{N_{1}\times N_{2}}$,
$\mathbf{B}_{1}\in\mathbb{C}^{N_{1}\times N_{2}}$, $c_{0}$, and
$c_{1}$ are generated in a manner analogous to Section \ref{subsec:Standard-Single-Equality-QCQP}. 

The performance comparison results are elaborated in Figure \ref{fig:main_comp_4_opt_time}.
Figure \ref{fig:Constraint-satisfaction-4} shows that validity ratio
of fmincon decreases with the problem size $N_{1}$ and eventually
falls below $50\%$. SDR generates valid solutions with a probability
of over $90\%$ but the problem size is no larger than 11. The solutions
of SD-BS are all valid across the entire range of problem sizes. As
for optimality gap, which is shown in Figure \ref{fig:Optimality-gap-4},
the proposed SD-BS is the smallest (down to $10^{-10}$) while the
other benchmarks are far from optimal and the average gap levels are
more than $10$. This could be because the SDR method reformulates
the original QCQP into a vectorized form with Kronecker-product patterns
which off-the-shelf solvers fail to exploit. With regard to computation
time, the proposed algorithm can be up to four orders of magnitude
faster than the benchmarks. The scalability range of SD-BS is at least
ten times that of the benchmarks. 

\section{Conclusion \label{sec:Conclusions}}

In this paper, we have studied a class of nonconvex single-equality-constraint
QCQP problems. Traditionally, this type of QCQP is solved via SDR
with high computational complexity. To improve algorithmic performance,
we have developed a non-iterative globally optimal algorithm named
SD-BS. The SD-BS algorithm follows a two-stage process: the SD stage
decomposes the vector-format variable into separable scalars with
an affine mapping and the BS stage computes the root of a monotonic
function obtained from first- and second-order KKT conditions. Moreover,
we have explored several extensions for applicability enhancement,
including rank deficiency, indefiniteness, linear constraint augmentation,
and matrix-format optimization variable. Numerical simulations have
proved the superiority of the proposed algorithms on general occasions
in terms of constraint satisfaction, optimality gap, and computational
efficiency. In terms of scalability, SD-BS can handle problem sizes
at least ten times as large as those manageable by off-the-shelf SDP
solvers. Specially in complex-valued matrix-format QCQP, SD-BS achieves
a lower objective value as well as a higher computational speed than
the solver-based benchmarks conditioned on feasibility. 

{\small

\bibliographystyle{IEEEtranN}
\bibliography{IEEEabrv,LichengZhao}

}

\begin{appendices}

\section{\label{sec:Proof-of-Theorem eqv recast}Proof of Theorem \ref{thm:eqv recast}}
\begin{IEEEproof}
We undo $\mathbf{A}_{1}$ with $\mathbf{T}\mathbf{T}^{T}$ and work
on the constraint of QCQP \eqref{eq: QCQP eq real} step by step:
\begin{equation}
\mathbf{x}^{H}\mathbf{A}_{1}\mathbf{x}=\mathbf{x}^{T}\mathbf{T}\mathbf{T}^{T}\mathbf{x}=\left(\mathbf{T}^{T}\mathbf{x}\right)^{T}\left(\mathbf{T}^{T}\mathbf{x}\right),
\end{equation}
\begin{equation}
2\mathbf{b}_{1}^{T}\mathbf{x}=2\mathbf{b}_{1}^{T}\mathbf{T}^{-T}\mathbf{T}^{T}\mathbf{x}=2\left(\mathbf{T}^{-1}\mathbf{b}_{1}\right)^{T}\left(\mathbf{T}^{T}\mathbf{x}\right),
\end{equation}
and thus 
\begin{equation}
\begin{aligned} & \mathbf{x}^{T}\mathbf{A}_{1}\mathbf{x}+2\mathbf{b}_{1}^{T}\mathbf{x}+c_{1}\\
= & \left(\mathbf{T}^{T}\mathbf{x}\right)^{T}\left(\mathbf{T}^{T}\mathbf{x}\right)+2\left(\mathbf{T}^{-1}\mathbf{b}_{1}\right)^{T}\left(\mathbf{T}^{T}\mathbf{x}\right)\\
 & +\left(\mathbf{T}^{-1}\mathbf{b}_{1}\right)^{T}\left(\mathbf{T}^{-1}\mathbf{b}_{1}\right)-\left(\mathbf{T}^{-1}\mathbf{b}_{1}\right)^{T}\left(\mathbf{T}^{-1}\mathbf{b}_{1}\right)+c_{1}\\
= & \left(\mathbf{T}^{T}\mathbf{x}+\mathbf{T}^{-1}\mathbf{b}_{1}\right)^{T}\left(\mathbf{T}^{T}\mathbf{x}+\mathbf{T}^{-1}\mathbf{b}_{1}\right)\\
 & -\left[\left(\mathbf{T}^{-1}\mathbf{b}_{1}\right)^{T}\left(\mathbf{T}^{-1}\mathbf{b}_{1}\right)-c_{1}\right]\\
= & \mathbf{y}^{T}\mathbf{y}-c
\end{aligned}
\end{equation}
with $\mathbf{y}=\mathbf{T}^{T}\mathbf{x}+\mathbf{T}^{-1}\mathbf{b}_{1}$
and $c=\left(\mathbf{T}^{-1}\mathbf{b}_{1}\right)^{T}\left(\mathbf{T}^{-1}\mathbf{b}_{1}\right)-c_{1}$.
With an affine mapping from $\mathbf{x}$ to $\mathbf{y}$, we work
on the objective analogously:
\begin{equation}
\mathbf{x}^{T}\mathbf{A}_{0}\mathbf{x}=\mathbf{x}^{T}\mathbf{T}\boldsymbol{\Lambda}_{0}\mathbf{T}^{T}\mathbf{x}=\left(\mathbf{T}^{T}\mathbf{x}\right)^{T}\boldsymbol{\Lambda}_{0}\left(\mathbf{T}^{T}\mathbf{x}\right),
\end{equation}
\begin{equation}
2\mathbf{b}_{0}^{T}\mathbf{x}=2\mathbf{b}_{0}^{T}\mathbf{T}^{-T}\mathbf{T}^{T}\mathbf{x}=2\left(\mathbf{T}^{-1}\mathbf{b}_{0}\right)^{T}\left(\mathbf{T}^{T}\mathbf{x}\right),
\end{equation}
and finally
\begin{equation}
\begin{aligned} & \mathbf{x}^{T}\mathbf{A}_{0}\mathbf{x}+2\mathbf{b}_{0}^{T}\mathbf{x}+c_{0}\\
= & \left(\mathbf{T}^{T}\mathbf{x}+\mathbf{T}^{-1}\mathbf{b}_{1}-\mathbf{T}^{-1}\mathbf{b}_{1}\right)^{T}\boldsymbol{\Lambda}_{0}\left(\mathbf{T}^{T}\mathbf{x}+\mathbf{T}^{-1}\mathbf{b}_{1}-\mathbf{T}^{-1}\mathbf{b}_{1}\right)\\
 & +2\left(\mathbf{T}^{-1}\mathbf{b}_{0}\right)^{T}\left(\mathbf{T}^{T}\mathbf{x}+\mathbf{T}^{-1}\mathbf{b}_{1}-\mathbf{T}^{-1}\mathbf{b}_{1}\right)+c_{0}\\
= & \left(\mathbf{T}^{T}\mathbf{x}+\mathbf{T}^{-1}\mathbf{b}_{1}\right)^{T}\boldsymbol{\Lambda}_{0}\left(\mathbf{T}^{T}\mathbf{x}+\mathbf{T}^{-1}\mathbf{b}_{1}\right)\\
 & -2\left(\boldsymbol{\Lambda}_{0}\mathbf{T}^{-1}\mathbf{b}_{1}-\mathbf{T}^{-1}\mathbf{b}_{0}\right)^{T}\left(\mathbf{T}^{T}\mathbf{x}+\mathbf{T}^{-1}\mathbf{b}_{1}\right)\\
 & +\left(\mathbf{T}^{-1}\mathbf{b}_{1}\right)^{T}\boldsymbol{\Lambda}_{0}\left(\mathbf{T}^{-1}\mathbf{b}_{1}\right)-2\left(\mathbf{T}^{-1}\mathbf{b}_{0}\right)^{T}\left(\mathbf{T}^{-1}\mathbf{b}_{1}\right)+c_{0}\\
= & \mathbf{y}^{T}\boldsymbol{\Lambda}_{0}\mathbf{y}-2\mathbf{b}^{T}\mathbf{y}+\mathrm{const}
\end{aligned}
\end{equation}
with $\mathbf{b}=\boldsymbol{\Lambda}_{0}\mathbf{T}^{-1}\mathbf{b}_{1}-\mathbf{T}^{-1}\mathbf{b}_{0}$
and $\mathrm{const}$ is an optimization irrelevant constant. 
\end{IEEEproof}

\section{\label{sec:Proof-of-Lemma SD-def}Proof of Lemma \ref{lem:SD-def}}
\begin{IEEEproof}
For any $\mathbf{A}_{1}\succeq\mathbf{0}$ of rank $r<N$, there exists
an invertible square matrix $\mathbf{S}_{\mathrm{A},1}$ such that
$\mathbf{S}_{\mathrm{A},1}\left[\begin{array}{cc}
\mathbf{I} & \mathbf{0}\\
\mathbf{0} & \mathbf{0}
\end{array}\right]\mathbf{S}_{\mathrm{A},1}^{T}=\mathbf{A}_{1}$ and the size of $\mathbf{I}$ is $r\times r$. Then, the multiplication
result of $\mathbf{S}_{\mathrm{A},1}^{-1}\mathbf{A}_{0}\mathbf{S}_{\mathrm{A},1}^{-T}\left(\succeq\mathbf{0}\right)$
is partitioned as \begin{equation}
\label{eq:decomp block}
\mathbf{S}_{\mathrm{A},1}^{-1}\mathbf{A}_{0}\mathbf{S}_{\mathrm{A},1}^{-T}=\mathbf{B}_{\mathrm{A,0}}=
\begin{bNiceMatrix}[last-row,last-col]    
 \mathbf{B}_{11}   & \mathbf{B}_{12}  & r \\
 \mathbf{B}_{21}   & \mathbf{B}_{22}  & N-r. \\
 r & N-r      \\ 
\end{bNiceMatrix}
\end{equation}The right-hand side of \eqref{eq:decomp block} is further decomposed
as 
\begin{equation}
\left[\begin{array}{cc}
\mathbf{B}_{11} & \mathbf{B}_{12}\\
\mathbf{B}_{21} & \mathbf{B}_{22}
\end{array}\right]=\mathbf{S}_{\mathrm{A},2}\left[\begin{array}{cc}
\mathbf{B}_{11}-\mathbf{B}_{12}\mathbf{B}_{22}^{\dagger}\mathbf{B}_{21} & \mathbf{0}\\
\mathbf{0} & \mathbf{B}_{22}
\end{array}\right]\mathbf{S}_{\mathrm{A},2}^{T},
\end{equation}
where $\mathbf{S}_{\mathrm{A},2}=\left[\begin{array}{cc}
\mathbf{I} & \mathbf{B}_{12}\mathbf{B}_{22}^{\dagger}\\
\mathbf{0} & \mathbf{I}
\end{array}\right]$. This decomposition is due to \citep[Corollary 8.2.2]{bernstein2009matrix}:
$\mathcal{R}\left(\mathbf{B}_{21}\right)\subseteq\mathcal{R}\left(\mathbf{B}_{22}\right)$
for $\mathbf{B}_{\mathrm{A,0}}\succeq\mathbf{0}$. Thus, $\mathbf{B}_{21}=\mathbf{B}_{22}\mathbf{B}_{22}^{\dagger}\mathbf{B}_{21}$
and $\mathbf{B}_{12}=\mathbf{B}_{12}\mathbf{B}_{22}^{\dagger}\mathbf{B}_{22}$.
Now we perform EVD on $\mathbf{B}_{11}-\mathbf{B}_{12}\mathbf{B}_{22}^{\dagger}\mathbf{B}_{21}$
and $\mathbf{B}_{22}$: $\mathbf{B}_{11}-\mathbf{B}_{12}\mathbf{B}_{22}^{\dagger}\mathbf{B}_{21}=\mathbf{V}_{1}\mathrm{Diag}\left(\mathbf{d}_{1}\right)\mathbf{V}_{1}^{T}$
and $\mathbf{B}_{22}=\mathbf{V}_{2}\mathrm{Diag}\left(\mathbf{d}_{2}\right)\mathbf{V}_{2}^{T}$.
Denote $\mathbf{S}_{\mathrm{A},3}=\left[\begin{array}{cc}
\mathbf{V}_{1} & \mathbf{0}\\
\mathbf{0} & \mathbf{V}_{2}
\end{array}\right]$ and eventually, simultaneous diagonalization is conducted as
\begin{equation}
\begin{cases}
\begin{aligned} & \mathbf{A}_{0}=\mathbf{S}_{\mathrm{A},1}\mathbf{S}_{\mathrm{A},2}\mathbf{S}_{\mathrm{A},3}\left[\begin{array}{cc}
\mathrm{Diag}\left(\mathbf{d}_{1}\right) & \mathbf{0}\\
\mathbf{0} & \mathrm{Diag}\left(\mathbf{d}_{2}\right)
\end{array}\right]\mathbf{S}_{\mathrm{A},3}^{T}\cdot\\
 & \mathbf{S}_{\mathrm{A},2}^{T}\mathbf{S}_{\mathrm{A},1}^{T}=\mathbf{T}\boldsymbol{\Lambda}_{0}\mathbf{T}^{T}
\end{aligned}
\\
\begin{aligned} & \mathbf{A}_{1}=\mathbf{S}_{\mathrm{A},1}\left[\begin{array}{cc}
\mathbf{I} & \mathbf{0}\\
\mathbf{0} & \mathbf{0}
\end{array}\right]\mathbf{S}_{\mathrm{A},1}^{T}=\mathbf{S}_{\mathrm{A},1}\mathbf{S}_{\mathrm{A},2}\left[\begin{array}{cc}
\mathbf{I} & \mathbf{0}\\
\mathbf{0} & \mathbf{0}
\end{array}\right]\mathbf{S}_{\mathrm{A},2}^{T}\mathbf{S}_{\mathrm{A},1}^{T}=\\
 & \mathbf{S}_{\mathrm{A},1}\mathbf{S}_{\mathrm{A},2}\mathbf{S}_{\mathrm{A},3}\left[\begin{array}{cc}
\mathbf{I} & \mathbf{0}\\
\mathbf{0} & \mathbf{0}
\end{array}\right]\mathbf{S}_{\mathrm{A},3}^{T}\mathbf{S}_{\mathrm{A},2}^{T}\mathbf{S}_{\mathrm{A},1}^{T}=\mathbf{T}\left[\begin{array}{cc}
\mathbf{I} & \mathbf{0}\\
\mathbf{0} & \mathbf{0}
\end{array}\right]\mathbf{T}^{T},
\end{aligned}
\end{cases}
\end{equation}
where $\mathbf{T}$ and $\boldsymbol{\Lambda}_{0}$ are redefined
as $\mathbf{S}_{\mathrm{A},1}\mathbf{S}_{\mathrm{A},2}\mathbf{S}_{\mathrm{A},3}$
and $\left[\begin{array}{cc}
\mathrm{Diag}\left(\mathbf{d}_{1}\right) & \mathbf{0}\\
\mathbf{0} & \mathrm{Diag}\left(\mathbf{d}_{2}\right)
\end{array}\right]$, respectively. The above mentioned process involves three EVD operations,
as opposed to two in Lemma \ref{lem:SD-full}. 
\end{IEEEproof}

\section{\label{sec:Proof-of-Theorem KKT mat}Proof of Theorem \ref{thm:KKT mat}}
\begin{IEEEproof}
For $\mathbf{y}_{i}$ with fixed element amplitudes, the phases of
$\mathbf{y}_{i}$ components should align on those of $\mathbf{b}_{i}$
so as to minimize the objective, i.e., 
\begin{equation}
\arg\left(\mathbf{y}_{i}\right)=\arg\left(\mathbf{b}_{i}\right),\forall i.
\end{equation}
Consequently, the decoupled QCQP \eqref{eq:decoupled QCQP summation}
is further simplified as:
\begin{equation}
\begin{aligned} & \underset{\left\{ \left|\mathbf{y}_{i}\right|\right\} _{i=1}^{N_{2}}}{\mathsf{minimize}} &  & \sum_{i=1}^{N_{2}}\mathbf{a}^{T}\left|\mathbf{y}_{i}\right|^{2}-2\left|\mathbf{b}_{i}\right|^{T}\left|\mathbf{y}_{i}\right|\\
 & \mathsf{subject\thinspace to} &  & \sum_{i=1}^{N_{2}}\mathbf{1}^{T}\left|\mathbf{y}_{i}\right|^{2}-c=0.
\end{aligned}
\end{equation}
The first-order KKT conditions are displayed below:
\begin{equation}
\sum_{i=1}^{N_{2}}\mathbf{1}^{T}\left|\mathbf{y}_{i}\right|^{2}=c\label{eq:norm sum equality}
\end{equation}
and 
\begin{equation}
\left|\mathbf{y}_{i}\right|=\left|\mathbf{b}_{i}\right|/\left(\mathbf{a}+\lambda\mathbf{1}\right).\label{eq:KKT y_i}
\end{equation}
Plugging \eqref{eq:decoupled QCQP summation} into \eqref{eq:decoupled QCQP summation},
we get a nonlinear equation $f\left(\lambda\right)=c$ where
\begin{equation}
\begin{aligned}f\left(\lambda\right)= & \sum_{i=1}^{N_{2}}\mathbf{1}^{T}\left[\left|\mathbf{b}_{i}\right|^{2}/\left(\mathbf{a}+\lambda\mathbf{1}\right)^{2}\right]\\
= & \sum_{i=1}^{N_{2}}\sum_{j=1}^{N_{1}}\frac{\left|b_{i,j}\right|^{2}}{\left(a_{j}+\lambda\right)^{2}}=\sum_{j=1}^{N_{1}}\frac{\sum_{i=1}^{N_{2}}\left|b_{i,j}\right|^{2}}{\left(a_{j}+\lambda\right)^{2}}.
\end{aligned}
\end{equation}
\end{IEEEproof}
\end{appendices}
\end{document}